\documentclass[12pt, uplatex]{amsart}
\usepackage{amsmath,amsfonts,amsthm,amscd,amssymb,mathrsfs,amssymb,bm,pb-diagram, here, xcolor}
\usepackage{setspace}
\usepackage[all,cmtip]{xy}

\usepackage{blkarray}
\usepackage{multirow}
\usepackage{mathtools}
\usepackage{graphicx}

\usepackage{mathrsfs, arydshln}
\usepackage{graphicx}
\usepackage[all]{xy}
\newtheorem{theorem}{Theorem}

\newtheorem{proposition}[theorem]{Proposition}
\newtheorem{lemma}[theorem]{Lemma}
\newtheorem{definition}[theorem]{Definition}

\newtheorem{corollary}[theorem]{Corollary}

\newtheorem{remark}[theorem]{Remark}
\newcommand{\aaa}{\alpha}

\newcommand{\CCC}{\Gamma}
\newcommand{\ddd}{\delta}
\newcommand{\DDD}{\Delta}

\newcommand{\lmd}{\lambda}
\newcommand{\Lmd}{\Lambda}

\newcommand{\CP}{\mathbb{CP}}
\newcommand{\PP}{\mathbb{P}}
\newcommand{\CC}{\mathbb{C}}

\newcommand{\RR}{\mathbb{R}}

\newcommand{\RP}{\mathbb{RP}}

\newcommand{\qdr}{\CP_1\times\CP_1}

\renewcommand{\hat}{\widehat}
\renewcommand{\tilde}{\widetilde}

\newcommand{\mf}{\mathfrak}

\newcommand{\ol}{\overline}

\newcommand{\lras}{\,\longrightarrow\,}

\newcommand{\set}{\,|\,}

\newcommand{\proofend}{\hfill$\square$}

\newcommand{\inv}{^{-1}}

\newcommand{\Aut}{{\rm{Aut}}}

\newcommand{\ms}{\mathscr}
\newcommand{\minus}{\backslash}
\newcommand{\ptl}{\partial}

\newcommand{\qandq}{\quad{\text{and}}\quad}

\newcommand{\reg}{{\rm{reg}}}
\newcommand{\pr}{{\rm{pr}}}

\newcommand{\us}{^{\sigma}}

\DeclareFontFamily{U}{mathx}{}
\DeclareFontShape{U}{mathx}{m}{n}{<-> mathx10}{}
\DeclareSymbolFont{mathx}{U}{mathx}{m}{n}
\DeclareMathAccent{\widehat}{0}{mathx}{"70}
\DeclareMathAccent{\widecheck}{0}{mathx}{"71}

\renewcommand{\hat}{\widehat}
\renewcommand{\tilde}{\widetilde}
\renewcommand{\check}{\widecheck}

\numberwithin{equation}{section}
\numberwithin{theorem}{section}

\begin{document}
\bibliographystyle{alpha} 
\title[]
{The Einstein-Weyl spaces associated to Segre quartic surfaces
}
\author{Nobuhiro Honda}
\address{Department of Mathematics,  Institute of Science Tokyo, 2-12-1, O-okayama, Meguro, 152-8551, JAPAN}
\email{honda@math.titech.ac.jp}
\author{Fuminori Nakata}
\address{Faculty of Human Development and Culture,
Fukushima University, JAPAN}
\email{fnakata@educ.fukushima-u.ac.jp}

\thanks{The first author was partially supported by JSPS KAKENHI Grant 16H03932.
}
\begin{abstract}
We find explicit examples of compact minitwistor spaces of genus one, whose Einstein-Weyl spaces have a connected component that is diffeomorphic to the de Sitter space.
The induced Einstein-Weyl structure on it is Lorenzian, real-analytic, whose spacelike geodesics are all closed and simple.
The identity component of the automorphism group of the Einstein-Weyl structure is the circle and therefore the structure is not isomorphic to the standard de Sitter structure. We show that these Einstein-Weyl structures deform as the Segre surfaces deform and converge to the standard de Sitter structure.

The minitwistor spaces we study are the so-called Segre quartic surfaces. They have a real pair of nodes, which play a crucial role in proving the above results. These singularities also allow us to construct explicit examples of non-compact complex surfaces that do not admit any compactification.
\end{abstract}
\maketitle

\setcounter{tocdepth}{1}
\section{Introduction}\label{s:I}
 
In this paper, we aim to uncover the intriguing differential geometry behind certain algebraic surfaces with a simple structure.

In \cite{Hi82}, Hitchin established a Penrose-type correspondence between certain complex surfaces known as the {\em minitwistor spaces} and 3-manifolds equipped with the Einstein-Weyl structure (EW structure for short).
The former are complex surfaces characterized by a family of smooth rational curves with self-intersection number two on them.
The EW structure, on the other hand, consists of a conformal structure and an affine connection that satisfy both a compatibility condition (known as the Weyl condition) and an Einstein equation.
The rational curves on a minitwistor space are referred to as {\em minitwistor lines}, and the space of these curves forms an EW 3-manifold.

This correspondence is local in nature, with its global version presented in \cite{HN11}, at least for the direction from minitwistor space to EW space.
In this generalization, minitwistor lines are still rational curves but may have several nodes.
The relationship between the number of nodes and the self-intersection number given there ensures that the space of minitwistor lines remains 3-dimensional.
Specifically, if $g$ denotes the number of nodes, the complete linear system generated by minitwistor lines induces a birational embedding of the minitwistor space into 
$\CP_{g+3}$ as a surface of degree $(2g+2)$ \cite{HN11}. 
Consequently, minitwistor lines on a minitwistor space are singular hyperplane sections of the surface under this embedding. This links the study of minitwistor lines to significant topics in algebraic geometry, such as the projective dual varieties and Severi varieties. We refer to the number 
$g$ of the nodes of a minitwistor line as the {\em genus} of the minitwistor space, because, under the embedding, a generic hyperplane section of the minitwistor space is a smooth curve of genus $g$.

In \cite{Hon22}, the first author showed that compact minitwistor spaces of {\em genus one} are precisely the Segre quartic surfaces, each of which is a complete intersection of two quadrics in $\CP_4$. 
These surfaces are referred to as {\em Segre surfaces} for short. 
They are classified into 16 types based on the normal forms of their defining equations \cite{Dol}, with one of these types consisting of smooth surfaces.

In this paper, we focus on a different type of Segre surfaces that have singularities. We introduce a real structure on these surfaces and investigate the associated real EW spaces.
These surfaces are obtained as a dimensional reduction of the (3-dimensional) twistor spaces of Joyce's self-dual metrics with torus action \cite{J95}, including the real structures.
Hence, the EW spaces associated to the surfaces have a connected component with a definite conformal structure that is the quotient of the self-dual structure. However, our interest in this paper is other connected component with Lorenzian EW structure.
We show that there exists a component that has the following remarkable properties.

\begin{theorem}\label{t:1}
There exists a 1-dimensional family of Segre surfaces with real structures for which the associated real EW spaces have a connected component with the following properties:
\begin{itemize}
\item The EW space is diffeomorphic to the de Sitter space 
$\mathbb S^2\times I$, where $I$ is an open interval.
\item The EW structure is Lorenzian and real analytic.
\item The EW structure is not isomorphic to the de Sitter structure but is spacelike Zoll (similarly to the de Sitter structure).
\item The identity component of the automorphism group of the EW structure is (not Lorenz group but) a circle.
\end{itemize}
\end{theorem}

For the third item, we call a Lorenzian EW space spacelike Zoll if all spacelike geodesics are closed, simple and moreover reduced in the sense that no geodesic is multiply covered by nearby geodesics.
The Segre surfaces and the real structures in the theorem are presented in explicit form in Section \ref{s:S}.

Our Segre surfaces depend on one real parameter and deform as complex surfaces when the parameter varies. However, it is not immediately clear whether the associated EW structure deforms in tandem with the Segre surface, since the global correspondence from EW space to minitwistor space has not been established.
We prove (without relying on this correspondence) that the associated EW structure does indeed deform as the Segre surface deforms. Furthermore, we clarify the behavior of the EW structures in the limit within the moduli space:

\begin{theorem}\label{t:2}
The EW spaces described in the previous theorem form a 1-dimensional moduli space, with the de Sitter structure appearing as a limit.
\end{theorem}

In particular, the de Sitter structure is deformable as an EW space, preserving not only Zoll property but also real analyticity. The deformability of the de Sitter structure as EW space has been known by \cite{LM09} and \cite{Nkt09}. In [13], the second author established such deformations satisfy the same Zoll property, but they only have low regularity in general. In \cite{LM09}, LeBrun and Mason established a one-to-one correspondence between orientation-reversing diffeomorphisms from $\CP_1$ to $\CP_1$ and EW spaces that satisfy some geometric conditions. In this correspondence, the complex conjugation map on $\CP_1$ corresponds to the de Sitter structure. It is not clear that the EW spaces constructed in Theorem 1.2 satisfy the geometric conditions in their theorem.

In the proof of these results, one of the key steps is to introduce a specific modification operation to an open subset of the Segre surface, which transforms any minitwistor lines in the open subset to minitwistor lines without nodes.
Hence, the modification transforms an open subset of the Segre surfaces into a minitwistor space in the original sense.
It turns out that the resulting minitwistor spaces do not admit a compactification.
In particular, we obtain the following result, which would be of independent interest:

\begin{theorem}\label{t:3} 
We can explicitly construct non-compact smooth complex surfaces that admit no compactification as a complex surface, from the Segre quartic surfaces.
\end{theorem}

As mentioned above, these complex surfaces have a structure of minitwistor space in the original sense. Non-compactifiability is shown by proving that the smooth minitwistor lines in the surface have $\ms O(2)$ as their normal bundles but they cannot be moved in a linear system.

\subsection{Outline of the paper and the proof of the theorems.}
Throughout this paper, the Segre surface under consideration is denoted by $S$. The defining equation of $S$ includes one real parameter $\kappa$ which satisfies $0<\kappa<1$
and $S$ always has two nodes as its only singularities.
In Section \ref{s:S}, after precisely defining $S$ in $\CP_4$ and the real structure $\sigma$ on it, we investigate the basic structures of $S$, such as a $\CC^*$-action, (rational) conic bundle structure,  double covering structures over, three smooth quadrics, and a double covering structure over 
the cone of a conic. In particular, the real locus of $S$ consists of two smooth spheres $S\us_1$ and $S\us_2$ which are `symmetric' to each other, and which are $\mathbb S^1$-invariant, where $\mathbb S^1\subset\CC^*$. 
These two spheres play a critical role throughout this paper.
At the end of this section, we explain why we focus on the present Segre surfaces and the real structure on them.

In Section \ref{s:rmtl}, among other things, we show that the space $W$ of real (i.e.\,$\sigma$-invariant) hyperplanes in $\CP_4$ that are tangent to $S$ at a point on the sphere $S\us_1$ and which intersect another sphere $S\us_2$ in a circle is diffeomorphic to a 3-manifold $\mathbb S^2\times I$.
This space $W$ has an $\mathbb S^1$-action induced from $S$, which is the product of the standard rotation on $\mathbb S^2$ and the trivial action on $I$. 
As a result, the $\mathbb S^1$-fixed locus in $W$ consists of two intervals.
We denote by $W_0$ the complement of these intervals in $W$.
It turns out that the hyperplane sections parameterized by $W_0$ are minitwistor lines, having a node at a point on $S\us_1$.
This implies that $W_0$ admits a Lorenzian EW structure (Proposition \ref{prop:W0}).
In contrast, all hyperplane sections parametrized by the two intervals are reducible curves consisting of one irreducible real conic and a real pair of two lines, which look like a fan; see Figure \ref{fig:fan}.
We call these reducible curves {\em irregular minitwistor lines}
(Definition \ref{d:imtl}).

For this reason, at first sight, it might be difficult to expect that the EW structure on $W_0$ extends smoothly to $W$.
In Sections \ref{ss:modify} and \ref{ss:extension} we prove that this is indeed the case.
In Section \ref{ss:modify}, we select a real connected open subset $U\subset S$ and introduce a series of operations to it, yielding a new open smooth complex surface $\check U$.
In Section \ref{ss:extension}, we show that nodal minitwistor lines and the irregular minitwistor lines contained in $U$ are transformed into {\em smooth} minitwistor lines in $\check U$, which means that $\check U$ is a minitwistor space in the original sense \cite{Hi82, JT85}. 
This implies the desired extension property of the EW structure (Theorem \ref{t:WEW}). 
Next, by investigating the dimension of the complete linear system generated by the (smooth) minitwistor lines in $\check U$, we provide a proof of Theorem \ref{t:3}. 
In the proof, a key role is played by pieces of curves that arise when we perform a transcendental operation in constructing $\check U$.
The fact that the present EW structure on $W=\mathbb S^2\times I$ 
is not isomorphic to the de Sitter structure readily follows from this (Corollary \ref{c:ndS}).
Finally, in Section \ref{ss:slg}, using the results obtained in Section \ref{s:rmtl}, we prove that our EW spaces are spacelike Zoll
and also the real analyticity of the EW structures
(Theorem \ref{thm:Zoll}).

In Section \ref{s:mad}, we prove Theorem \ref{t:2}.
In Sections \ref{ss:disk} and \ref{ss:inv}, we investigate the family of disks that are halves of irreducible minitwistor lines divided by the real circles. We prove that the complement of two specific straight lines in the Segre surface admits a double fibration \eqref{d:double3}.
The fibers of one fibration are disks, and the fibers of another fibration are either real null surfaces or timelike geodesics, depending on whether they are over a real or a non-real point of the minitwistor space respectively.
Next, in Section \ref{ss:var2}, using this double fibration and a slight generalization of the modification in Section \ref{ss:modify}, we prove that the EW structures obtained from two Segre surfaces $S$ and $S'$ are isomorphic if and only if they are biholomorphic preserving the real structure (Theorem \ref{t:var}), which proves the first part of Theorem \ref{t:2}.
Finally, in Section \ref{ss:dS}, we show that the present EW spaces converge to the de Sitter space.
Broadly speaking, this is shown by first applying a linear coordinate change involving the parameter $\kappa\in (0,1)$ to change the defining equations of $S$ and then letting $\kappa$ be zero.
In the limit, $S$ will split into a cone $T_1$ over a conic and a smooth quadric $T_2$ which is the minitwistor space of the de Sitter space. Every minitwistor line in $S$ will split into a smooth conic in $T_2$ and two generating lines in $T_1$, where the real node of the minitwistor line converges to the vertex of the cone $T_1$.
Utilizing the double fibration obtained in Sections \ref{ss:disk} and \ref{ss:inv}, we show the convergence to the de Sitter space.

\section{Some Segre surfaces equipped with a real structure}
\label{s:S}
Consider a complete intersection of two quadrics in $\CP_4$
defined by the following equations
\begin{align}\label{Segre1}
X_0 X_1 - \kappa X_2^2 +  X_3^2 = 
X_2^2 -  X_3^2 -  X_4^2 = 0.
\end{align}
Here, $X_i$ are homogeneous coordinates on $\CP_4$ and $\kappa$ is any real number satisfying $0<\kappa<1$.
Throughout this article, we denote this surface by $S$.
These are among a special type of so-called Segre quartic surfaces.
For simplicity, we refer to them as Segre surfaces in this paper.
The surface $S$ has two ordinary nodes 
\begin{align}\label{e_0e_1}
e_0:=(1,0,0,0,0)\qandq e_1:=(0,1,0,0,0)
\end{align} 
and these are all singularities of $S$. 
Further, $S$ is invariant under the real structure
\begin{align}\label{rs1}
\sigma:(X_0,X_1,X_2,X_3,X_4)\longmapsto
\big(\ol X_1,\ol X_0,\ol X_2,\ol X_3,\ol X_4\big).
\end{align}
Then $\sigma(e_0) = e_1$.
The surface $(S,\sigma)$, which depends on a real number $\kappa\in(0,1)$, is the object we shall consider throughout this paper. These surfaces are anti-canonical models of quartic del Pezzo surfaces.

The $\CC^*$-action on $\CP_4$ defined by 
\begin{align}\label{C*action1}
(X_0,X_1,X_2,X_3,X_4)
\longmapsto
\big(tX_0,t\inv X_1,X_2,X_3,X_4\big),\quad t\in\CC^*
\end{align}
preserves $S$, constitute the identity component of the holomorphic automorphism group of $S$.
The map \eqref{C*action1} commutes with the real structure \eqref{rs1} if and only if $|t|=1$. So $\Aut_0(S,\sigma)=\mathbb S^1$.

Let $f:\CP_4\lras \CP_2$ be the projection 
from the line $\ol{e_0e_1}=\{X_2 = X_3 = X_4 = 0\}$.
This preserves the real structure $\sigma$.
From the second equation of \eqref{Segre1}, 
the image of $S$ under $f$ is the conic
$$
\Lmd:= \big\{(X_2,X_3,X_4)\in\CP_2\set X_2^2 -  X_3^2 -  X_4^2 = 0\big\}.
$$
We use the same symbol $f$ to mean 
$f|_S$, so that we have a dominant rational map
\begin{align}\label{pi}
f:S\lras \Lmd.
\end{align}
This map $f$ realizes $S$ as a rational conic bundle 
over $\Lmd$, with the singular points $e_0,e_1$ as its indeterminacy locus.
From the equations \eqref{Segre1},
the fiber conic $f\inv(\lmd)$ over a point 
$\lmd\in\Lmd$ degenerates into two lines
iff $\lmd$ is one of the following four points
\begin{multline}
\lmd_1:=\big(1, \sqrt{\kappa},\sqrt{1-\kappa}\big),\,\,
\lmd_2:=\big(1, -\sqrt{\kappa},\sqrt{1-\kappa}\big),\\
\lmd_3:=\big(1, -\sqrt{\kappa},-\sqrt{1-\kappa}\big),\,\,
\lmd_4:=\big(1, \sqrt{\kappa},-\sqrt{1-\kappa}\big).
\label{disc1}
\end{multline}
We call these the {\em discriminant points} of $f$.
All these are real
for the real structure
$(X_2,X_3,X_4)\longmapsto(\ol X_2,\ol X_3,\ol X_4)$
induced from $\sigma$. 
In real non-homogeneous coordinates
$(X_3/X_2,X_4/X_2)$, the real locus $\Lmd^{\sigma}$ is 
the unit circle and the four points are placed as in Figure \ref{fig:circle1}.
Each of the four fibers $f\inv(\lmd_i)\subset S$ consists of two lines. We write $f\inv(\lmd_i) = l_i + \ol l_i$, where
we make distinction of the two lines by supposing $e_0\in l_i$ and $e_1\in \ol l_i$.
These eight are all lines on $S$. 
We also put 
\begin{align}\label{p_i}
p_i:=l_i\cap \ol l_i,\quad 1\le i\le 4.
\end{align}
These are real points and $f(p_i) = \lmd_i$.
The 6 points $e_0,e_1,p_1,\dots,p_4$ are all $\CC^*$-fixed points on $S$.
We may use the (rational) conic bundle structure on $S$ to determine the real locus:

\begin{figure}
\includegraphics{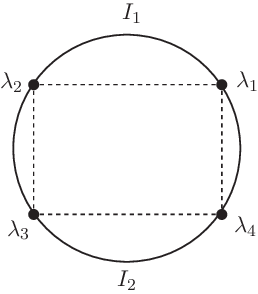}
\caption{
The real locus $\Lmd^{\sigma}$ 
}
\label{fig:circle1}
\end{figure}

\begin{proposition}\label{prop:Ssigma}
The real locus $S^{\sigma}$ of $S$ 
consists of two disjoint 2-spheres, smoothly embedded in $S\minus\{e_0,e_1\}$,
and they are lying over the arcs in the circle $\Lmd\us$ whose endpoints are $\{\lmd_1,\lmd_2\}$ and $\{\lmd_3,\lmd_4\}$ respectively.
\end{proposition}

\proof 
By a coordinate change $X_0= Y_0+\sqrt{-1}Y_1$ and $X_1= Y_0-\sqrt{-1}Y_1$, the real structure $\sigma$ becomes just the complex conjugation and the equations of $S$ become 
$Y_0^2 + Y_1^2 -\kappa X_2^2 + X_3^2 = X_2^2 - X_3^2 - X_4^2=0$.
From this, $S\us\cap \{X_2=0\}=\emptyset$ and
in the real affine coordinates $(y_0,y_1,x_3,x_4) = 
(Y_0,Y_1,X_3,X_4)/X_2$,
$S\us$ is defined by 
$$
y_0^2 + y_1^2 + x_3^2 = \kappa,\quad x_3^2 + x_4^2 = 1.
$$
This consists of two disjoint spheres 
$
\{y_0^2 + y_1^2 + x_3^2 = \kappa,\, x_4>0\}
$ and 
$\{y_0^2 + y_1^2 + x_3^2 = \kappa,\, x_4<0\}.
$
In the coordinates, the circle $\Lmd\us$ is defined by $x_3^2 + x_4^2 = 1$ and the latter assertion follows immediately from this.
\proofend

\begin{definition}\label{def:Ssigma}
{\em
We denote $I_1:=(\lmd_1,\lmd_2)$ and $I_2:= (\lmd_3,\lmd_4)$
for the open arcs in the circle $\Lmd^{\sigma}$.
(See Figure \ref{fig:circle1}.)
We denote $S^{\sigma}_1$ and $S^{\sigma}_2$ for
the connected components of $S^{\sigma}$,
 which are 2-spheres lying over the 
the closures of $I_1$ and $I_2$ in $\Lmd^{\sigma}$ respectively.
\proofend
}
\end{definition}

The spheres $S^{\sigma}_1$ and $S^{\sigma}_2$ play important roles throughout this article.
They are $\mathbb S^1$-invariant, where $\mathbb S^1\subset\CC^*$.
By definition, $p_1,p_2\in S^{\sigma}_1$ and  $p_3,p_4\in
S^{\sigma}_2$, and the $\mathbb S^1$-actions on these spheres are rotations around these points. 

To each $j=2,3,4$, let $\pi_j:\CP_4\lras\CP_3$ be the projection
which drops the coordinate $X_j$.
This induces by restriction a finite degree-two mapping onto a smooth quadric $Q_j:=\pi_j(S)\subset\CP_3$.
So $S$ has three double covering structures over smooth quadrics.
Let $B_j\subset Q_j$ be the branch curve.
Each of these consists of two irreducible $(1,1)$-curves intersecting at two points, and these points are the images of the singularities $e_0$ and $e_1$.
The real structure and the $\CC^*$-action on $S$ descend on all the three quadrics and the branch curves,
and generic orbits of the induced $\CC^*$-action on $Q_j$ are 
smooth $(1,1)$-curves through the two points $\pi_j(e_0)$ and $\pi_j(e_1)$.
The $\CC^*$-action on $Q_j$ has two other fixed points
and they are real points. These are the images of the four points
$p_1,\dots,p_4$.
If $\tau_j:S\lras S$ ($j=2,3,4$) is the holomorphic involution which exchanges $X_j$ and $-X_j$, then
the map $\pi_j:S\lras Q_j$ can be regarded as the quotient map by $\tau_j$.
From this, we can detect which pair of points among $p_1,\dots,p_4$ are lying over the same point under $\pi_j$.

\begin{proposition}\label{p:Qrs}
The induced real structure on $Q_j\simeq\qdr$ ($j=2,3,4$) is of the form $(z,w)\longmapsto (\ol w,\ol z)$
for some affine coordinates $(z,w)$ on $\qdr$. In particular,
$Q_j^{\sigma}\simeq\mathbb S^2$.
\end{proposition}

We omit a proof of this proposition because it can be readily shown from the explicit description of the objects.
Pulling back by $\pi_j:S\lras Q_j$ the pencils of $(1,0)$-lines and $(0,1)$-lines on the quadrics $Q_j\simeq\qdr$, we obtain six pencils on $S$. 
These pencils are mutually distinct and generic members of these pencils are irreducible conics.
The $\CC^*$-action \eqref{C*action1} on $S$ preserves each of these pencils and the pencils have exactly two $\CC^*$-invariant members, all of which consist of two lines.
These invariant members can be explicitly determined using the four points $p_1,\dots, p_4$ and the above involutions $\tau_j$.
From Proposition \ref{p:Qrs}, for each $j=2,3,4$, the real structure on $S$ exchanges the two pencils.
These pencils will be used in Section \ref{ss:disk}.

\begin{proposition}
\label{prop:conics}
None of these six pencils on $S$ are linearly equivalent to each other.
Any irreducible conic on $S$ is either a member of these pencils or a fiber of the (rational) conic bundle map $f:S\lras\Lmd\simeq\CP^1$.
\end{proposition}

This can be shown using the following realization of $S$ as a birational transformation of $\qdr$.
Take two points on the same $(1,0)$-curve on $\qdr$
and also two points on another $(1,0)$-curve, requiring that no two points among these four points
are on the same $(0,1)$-curve.
Next, let $\tilde S\lras \qdr$ be the blowing up at these four points.
The strict transform of the two $(1,0)$-curves are $(-2)$-curves on $\tilde S$, so they can be blown down to ordinary nodes.
The resulting surface is biholomorphic to one of the present Segre surfaces $S$.
The $\CC^*$-action on $S$ is obtained as the transform of the $\CC^*$-action on $\qdr$ whose orbits closures are $(0,1)$-curves.
The eight lines on $S$ are the four exceptional curves and the strict transforms of the four $(0,1)$-curves passing through the blownup points.
The freedom of choosing the four points corresponds to the choice of the parameter $\kappa$ in the equation \eqref{Segre1} of $S$.
On the other hand, the real structure $\sigma$ given in \eqref{rs1} can never be obtained as a transformation of that on $\qdr$
as it switches $l_i$ and $\ol l_i$.

\begin{proposition}\label{prop:2sing}
Any hyperplane whose section of $S$ has at least two singularities that are different from $e_0$ and $e_1$ has to be of one of the following forms:
\begin{itemize}
\item[(i)] $\pi_j\inv(T_q Q_j)$, $j=2,3,4$ and $q\not\in B_j$,
\item[(ii)] $f\inv(l)$, where $l$ is a tangent line of the conic $\Lmd$,
\item[(iii)] $f\inv(l)$, where $l$ is the line $\ol{\lmd_i\lmd_j}$
for some $1\le i<j\le 4$.
\end{itemize}
\end{proposition}

Again we omit a proof
as it can be shown by a standard argument for a projective surface. 
This will be used when we classify minitwistor lines on $S$.

The following proposition is about the real structure induced on the quadrics $Q_j$.
\begin{proposition}\label{prop:3rs}
The real structures on $Q_j$ satisfy the following properties.
\begin{itemize}
\item
The real structure on $Q_2$ exchanges the two components of
the branch curve $B_2$,
while those on $Q_3$ and $Q_4$ preserve the two components
of $B_3$ and $B_4$ respectively.
Further, we have 
$$
B_2^{\sigma} = B_4^{\sigma} = \emptyset
\qandq
B_3^{\sigma} = \mathbb S^1\sqcup \mathbb S^1,
$$
where the last two circles are contained in the distinct components
of $B_3$.
\item
When $j=2,4$, 
the restriction of the double covering $\pi_j:S\lras Q_j$ to the real locus 
$S^{\sigma} = S^{\sigma}_1\cup S^{\sigma}_2$ is identified with the trivial double cover over $Q_j^{\sigma}\simeq \mathbb S^2$.
When $j=3$, the same restriction is not surjective over $Q^{\sigma}_3$, and 
the images $\pi_3(S^{\sigma}_1)$ and 
$\pi_3(S^{\sigma}_2)$ are
mutually disjoint closed disks in $Q^{\sigma}_3\simeq \mathbb S^2$.
The two spheres $S^{\sigma}_1$ and $S^{\sigma}_2$ are mapped
to these two disks as a degree-two map
whose branch locus is the boundaries of the disks.
\end{itemize}
\end{proposition}

Again we omit a proof because it can be seen from the explicit description.

Next, let $P\in\CP_4$ be any {\em real} point on the line $\ol{e_0e_1}$.
Since $S\cap \ol{e_0e_1} = \{e_0,e_1\}$ and $e_1 = \sigma(e_0)$, $P\not\in S$.
Let $\Pi:S\lras\CP_3$ be the restriction of the projection $\CP_4\lras\CP_3$ from $P$ to $S$. This is holomorphic and the image $\Pi(S)\subset\CP_3$ is the cone over the conic $\Lmd$.
We denote $C(\Lmd)$ for this cone, so we have a holomorphic map $\Pi:S\lras C(\Lmd)$. It is easy to see that $\Pi$ realizes $S$ as a double covering over $C(\Lmd)$ and its branch divisor is the intersection of $S$ with a quadric. 
If $\Sigma\subset C(\Lmd)$ is the branch divisor, then $\Sigma$ is a smooth elliptic curve and it does not pass through the vertex of $C(\Lmd)$.
The nodes $e_0$ and $e_1$ of $S$ are over the vertex of the cone.
Since $P$ is supposed to be real, $\Pi$ preserves the real structure and $\Sigma$ is real.
However, since $P$ is not fixed by the $\CC^*$-action on $\CP_4$ of \eqref{C*action1}, unlike $\pi_j$ ($j=2,3,4$), the $\CC^*$-action on $S$ does not descend on $C(\Lmd)$.
On the other hand, since $\ol{e_0e_1}\us\simeq \mathbb S^1$ and this is an orbit of the $\mathbb S^1$-action, the structure of the double covering $\Pi$ does not depend on the choice of the point $P$.
This double covering structure $\Pi$ will be used in Section \ref{ss:var2}.

\begin{remark}\label{rmk:Joyce}
{\em
Because there are many types of Segre surfaces, and they do not appear to have a natural real structure, we explain why we focus on the one considered here.
Jones and Tod \cite{JT85} demonstrated that a minitwistor space can be obtained as a 1-dimensional reduction from a 3-dimensional twistor space. Applying this to the twistor spaces of Joyce's self-dual metrics with torus action \cite{J95}, the first author obtained several compact minitwistor spaces as explicit projective rational surfaces \cite{Hon10}.
A smooth quadric is the simplest one obtained in this way, and the present Segre surfaces are the second simplest.
The real structure \eqref{rs1} on these Segre surfaces is induced by the real structure on the 3-dimensional twistor spaces.
Moreover, the $\mathbb{S}^1$-action on $S$, obtained by setting $|t| = 1$ in \eqref{C*action1}, is the residual action of the torus of the metric.

The image of a generic real twistor line on the twistor space to the present surface, which is the minitwistor line in the surface, has exactly one node as its only singularity, and this node is the unique real point on the minitwistor line.
The absence of a real circle on the minitwistor line reflects the fact that Joyce's structure is (positive) definite.
However, the main focus of this paper is not these minitwistor lines but rather those whose real locus contains a real circle.
This implies that the EW structures on the 3-manifold will be indefinite.
Therefore, the EW structures considered here are {\em not} obtained from Joyce's self-dual metrics by dimensional reduction.
}
\end{remark}

\section{Real minitwistor lines on the Segre surfaces}
\label{s:rmtl}
According to \cite[Theorem 2.7]{Hon22}, any Segre surface is a minitwistor space of genus one.
Thus, by \cite[Theorem 2.10]{HN11}, the space of minitwistor lines on $S$ is a 3-dimensional complex manifold equipped with a complex EW structure.
If the Segre surface has a real structure, the space of real minitwistor lines becomes a 3-dimensional real EW manifold.
In this section, we investigate this space for the Segre surfaces with a real structure, as studied in the previous section.

Let $S \subset \mathbb{CP}_4$ be a Segre surface \eqref{Segre1}, and let $\sigma$ denote the real structure on $S$ given by \eqref{rs1}.
The minitwistor lines on $S$ are irreducible rational curves contained in the smooth locus $S_{\text{reg}}$ that have a self-intersection number of four and a single node as their only singularity.
These lines are obtained as hyperplane sections of $S$ that are tangent to $S$ at a smooth point (see \cite[Section 2.1]{Hon22}).
Note that the section of $S$ by such a hyperplane is not always a minitwistor line, as it can be a reducible curve or have a singularity that is not an ordinary node (e.g., a cusp) at the tangent point.
If $C$ is a {\em real} minitwistor line on $S$, then the node of $C$ must be a real point of $S$, due to the uniqueness of the node.
Thus, any real minitwistor line on $S$ takes the form $S \cap H$, where $H \subset \mathbb{CP}^4$ is a real (i.e., $\sigma$-invariant) hyperplane that is tangent to $S$ at a real (and hence smooth) point of $S$.

For any real point $p$ on $S$, let $(T_pS)^{*}$ denote the pencil formed by hyperplanes containing the tangent plane $T_pS$, and let $(T_pS)^{*\sigma}$ denote the real locus of $(T_pS)^{*}$, which is always a circle. We use $(T_pS)^*|_S$ to refer to the pencil on $S$ consisting of sections of $S$ by hyperplanes in $(T_pS)^{*}$, and similarly for $(T_pS)^{*\sigma}|_S$. The goal of this section is to determine, for an arbitrary real point $p$ on $S$, the structure of an arbitrary member of the real pencil $(T_pS)^{*\sigma}|_S$ and its real locus.

\subsection{Real hyperplane sections which are not minitwistor lines}\label{ss:nmtl}

In this subsection, for the above goal, to any real point $p\in S^{\sigma}$, we determine all members of the real pencil $(T_pS)^{*\sigma}|_S$ which are {\em not} minitwistor lines. If $p$ is one of the $\mathbb S^1$-fixed points $p_1,\dots,p_4$, then since these are intersection points of two lines on $S$ (see \eqref{p_i}), any member of $(T_pS)^{*\sigma}|_S$ has these lines as components. This case is of another significance and will be treated in Section \ref{ss:imtl}. In this section we mostly assume $p\not\in\{p_1,\dots,p_4\}$.

From Proposition \ref{prop:Ssigma},
the real locus $S\us$ consists of two spheres $S^{\sigma}_1$ and $S^{\sigma}_2$. 
The next proposition indicates that it is sufficient to consider only real minitwistor lines whose nodes lie on $S^{\sigma}_1$. This can be immediately shown from the definitions of $\tau_2$ and $\tau_4$ and \eqref{disc1}, so we omit the proof.

\begin{proposition}\label{p:symm}
The holomorphic involutions $\tau_2$ and $\tau_4$ on $S$
exchange $S^{\sigma}_1$ and $S^{\sigma}_2$,
and they satisfy $\tau_2(p_1) = p_3,\, \tau_2(p_2) = p_4,\,
\tau_4(p_1) = p_4,\, \tau_4(p_2) = p_3.$
\end{proposition}

Since $S$ is a surface of degree four in $\CP_4$, a hyperplane section is a quartic curve, and a generic one is a smooth elliptic curve.

\begin{proposition}\label{prop:irr1}
If $p\in S\us_1\minus\{p_1,p_2\}$, then 
a generic member of the pencils $(T_pS)^*|_S$ and $(T_pS)^{*\sigma}|_S$ is irreducible.
\end{proposition}

\proof
Since any line on  $S$ does not pass through $p$ under the assumption,  a member of 
the pencil $(T_pS)^*|_S$ can be reducible only when it consists of two irreducible conics.
Each of these two conics has to pass through the point $p$
since otherwise $S\cap H$ would be smooth at $p$.
But this cannot happen for a generic $H\in (T_pS)^*$,
because any smooth conic on $S$ which does not pass through singularities of $S$ has self-intersection number zero and
therefore it cannot fix any point.
\proofend

\medskip
We will soon show that the singularity at $p$ of a generic member of the pencil is a node.
Hence, a generic member of the pencil $(T_pS)^*|_S$ will be a (real) minitwistor line. 
On the other hand, using the projections $f:\CP_4\lras\CP_2$ and $\pi_j:\CP_4\lras\CP_3$ as in the previous section,
it is easy to find real hyperplane sections of $S$ which are not minitwistor lines:

\begin{definition}\label{d:ehpls}{\em
For any point $p\in S\us_1$ (including the case $p=p_1,p_2$),
we define four real hyperplanes by
\begin{align}\label{H_0}
H_0(p)&:= f\inv\big(T_{f(p)}\Lmd\big),\\
\label{H_j}
H_j(p)&:=\pi_j\inv(T_{q_j}Q_j),\quad q_j := \pi_j(p),
\quad j=2,3,4.
\end{align}
}
\end{definition}

The hyperplane $H_0(p)$ is tangent to $S$ along the fiber conic over the point $f(p)\in\Lmd$. So $H_0(p)\in (T_pS)^{*\sigma}$ and $S\cap H_0(p)$ is a double conic. 
If $j=2,3,4$, then since $Q_j\cap T_{q_j}Q_j\subset\CP_3$ consists of two lines, $S\cap H_j(p)$ is reducible and has $p$ as one of its singularity.
So $H_j(p)$ ($j=2,3,4$) belong to the pencil $(T_pS)^{*\sigma}$ and
$S\cap H_j(p)$ are not minitwistor lines.
Using Proposition \ref{prop:conics}, the involutions $\tau_j$, and also the above reducibility, 
we can readily show that these four hyperplanes are distinct.
Further, for the structure of the sections $S\cap H_j(p)$, we readily obtain from Proposition \ref{prop:3rs} the following.
Recall that $B_j$ is the branch divisor of $\pi_j:S\lras Q_j$.

\begin{lemma}\label{lemma:nonmtl1}
Suppose $p\in S\us_1\minus\{p_1,p_2\}$.
To any $j=2,3,4$, the real hyperplane section $ S \cap H_j(p)$
always consists of two irreducible conics,
and these are exchanged by $\sigma$.
Further, if $j=2,4$, then the two conics always intersect transversely at the two points
$p$ and $\tau_j(p)$.
The same conclusion holds in the case $j=3$ unless $\pi_3(p)\in B_3$, 
in which case the two conics are tangent to each other
at the point $p$. 
\end{lemma}

Note that, since $B_2^{\sigma} = B_4^{\sigma} = \emptyset$ by Proposition \ref{prop:3rs}, $\pi_j(p)\in B_j\us$ can happen only when $j=3$.
Note also that the tangential situation of the two conics
happens as a limit of the transversal situation,
by just approaching a generic point of $S^{\sigma}_1$ to 
a real point of the ramification curve of $\pi_3$.

Next, we find yet another member of the real pencil 
$ (T_pS)^{*\sigma}|_S$ which are not a minitwistor line.

\begin{lemma}\label{lemma:nonmtl2}
For any point $p\in S\us_1\minus\{p_1,p_2\}$, apart from the double conic $H_0(p)|_S$, there exists a unique member
of the pencil $(T_pS)^{*}|_S$
which has a non-nodal singularity at $p$, and the member is real.
This singularity is an ordinary cusp unless $\pi_3(p)\in B_3^{\sigma}$.
\end{lemma}

Note that by Lemma \ref{lemma:nonmtl1}, when $\pi_3(p)\in B_3^{\sigma}$,
the hyperplane section $S\cap H_3(p)$ has a tacnode at $p$.
This means that if $\pi_3(p)\in B_3^{\sigma}$, then the unique hyperplane in the lemma is exactly $H_3(p)$.

\medskip
\noindent {\em Proof of Lemma \ref{lemma:nonmtl2}.}
The former half of the assertion is in effect proved in \cite[Lemma 3.7 and Proposition 3.8]{cuspidal}, but we give an outline of the proof since 
in \cite{cuspidal} no real structure is taken into account, and also because later we will use a part of the proof.

Take non-homogeneous coordinates $(x,y,z,w)$ in $\CP_4$ centered at 
the point $p$,
such that the real structure $\sigma$ takes the form
$(x,y,z,w)\longmapsto (\ol x, \ol y, \ol z, \ol w)$.
Without loss of generality, we may suppose that, in a neighborhood of $p$, $S$ is defined by 
the equations $F=F(x,y)$ and $G=G(x,y)$,
where $F$ and $G$ are holomorphic functions 
that are real-valued when $x$ and $y$ are real.
Any hyperplane $H$ belonging to the pencil $(T_pS)^{*}$ is
defined by an equation of the form
\begin{multline}\label{tpx2}
\lmd\left\{z - F(0,0)
- 
F_x(0,0) x
- 
F_y(0,0) y\right\}\\
+
\mu
\left\{
w - G(0,0)
- 
G_x(0,0) x
- 
G_y(0,0) y
\right\}=0
\end{multline}
for some $(\lmd, \mu)\in \CC^2\minus\{(0,0)\}$,
and $H$ is real iff $\lmd$ and $\mu$ can be taken as real numbers.
Substituting $z = F(x,y)$ and $w = G(x,y)$ into \eqref{tpx2},
we obtain a local equation of the hyperplane section $S\cap H$ around $p$.
Let ${\bf H}={\bf H}(\lmd,\mu)$ be the Hessian of the last equation with respect to $x$ and $y$.
The function $\bf H$ is a real quadratic polynomial in $\lmd$ and $\mu$ unless it is identically zero.
(See \cite[(2.3)]{cuspidal} for the concrete form of ${\bf H}$.)

As mentioned in the proof of Proposition \ref{prop:irr1},
since $p\neq p_i\,(i=1,\dots, 4)$, there exists no line on $S$
which passes through $p$.
By \cite[Proposition 3.8]{cuspidal}, this means that the function ${\bf H}(\lmd,\mu)$ is not identically zero and does not have a double root as a function on $\CP_1$.
This implies that generic members of the pencil 
$(T_pS)^{*}|_S$ have ordinary nodes at $p$, and there exist
precisely two members whose singularities at $p$ are not ordinary nodes.
The double conic $H_0(p)|_S$ is one of the last two members of $(T_pS)^{*}|_S$.
Let $H_1(p)$ be the other member of $(T_pS)^{*}$.
Since the pencil $(T_pS)^{*}$ is real and $H_0(p)$ is real, $H_1(p)$ is also real. Thus, we have shown the uniqueness in the lemma.

Finally, we determine the singularity of $S\cap H_1(p)$ at $p$
when $\pi_3(p)\not\in B_3^{\sigma}$.
Suppose that 
the section $S\cap H_1(p)$ is reducible.
As in the proof of Proposition \ref{prop:irr1},
we can then write $S\cap H_1(p) = C_1 + C_2$,
where $C_1$ and $C_2$ are mutually distinct smooth conics.
If one of them, say $C_1$, is a fiber of the rational map
$f:S\to \Lmd$, then $e_0,e_1\in C_1$, which means 
$\ol{e_0e_1}\subset H_1(p)$.
But the section of $S$ by such a hyperplane can have a non-nodal singularity at $p$ only when it is exactly a double conic $H_0(p)|_S$, but it cannot happen since $H_1(p)\neq H_0(p)$.
Therefore, $S\cap H_1(p)$ is irreducible.
Then as the arithmetic genus of the hyperplane section of $S$ is one, 
the singularity of $S\cap H_1(p)$ at $p$ has to be an ordinary cusp.
This completes a proof of Lemma \ref{lemma:nonmtl2}.
\proofend

\medskip
From the lemma,  for any point $p\in S\us_1\minus\{p_1,p_2\}$,
we call the unique element of the real pencil
$(T_pS)^{*\sigma}|_S$ in the lemma
the {\em cuspidal member}, and write $H_1(p)$
for the real hyperplane which cuts out the cuspidal member.
As above, when $p$ satisfies 
$\pi_3(p)\in B_3^{\sigma}$, the singularity of the section $S\cap H_1(p)$ at $p$
is not a cusp but a tacnode.
But even in this case, we call $S\cap H_1(p)$ a cuspidal member,
because in the case $\pi_3(p)\in B_3^{\sigma}$, $S\cap H_3(p)$ can be considered as a limit of a generic cuspidal member.

\begin{lemma}\label{lemma:nonmtl3}
For any $p\in S\us_1\minus\{p_1,p_2\}$, the hyperplanes $H_i(p)$ ($0\le i\le 4$) 
are all members of the real pencil $(T_pS)^{*\sigma}$ which do not cut out minitwistor lines from $S$.
\end{lemma}

\proof
Suppose that $H\in (T_pS)^{*\sigma}$ is such that
$S\cap H$ is not a minitwistor line.
If $H\neq H_0(p)$ and $H\neq H_1(p)$, then
from the uniqueness as in Lemma \ref{lemma:nonmtl2}, the singularity of $S\cap H$ at $p$ is an ordinary node.
As $S\cap H$ is not a minitwistor line, this means that 
$S\cap H$ has another singularity.
From Proposition \ref{prop:2sing}, this means that $H$ is one of the three possibilities in the same proposition.
But among them, only (i) satisfies the present situation.
This implies $H= H_j(p)$ for some $j=2,3,4$.
\proofend

\subsection{Transition in the structure of the real locus of 
real hyperplane sections}\label{ss:var1}
In this subsection, still fixing any point $p\in S\us_1\minus\{p_1,p_2\}$, we investigate how the real locus of the 
hyperplane section $S\cap H$ varies as $H$ moves in the real pencil
$(T_pS)^{*\sigma}$.
Obviously, the topological structure of the real locus can change only when
$H$ passes through an $H$ for which $S\cap H$ is not a minitwistor line. Namely, by Lemma \ref{lemma:nonmtl3}, exactly when $H = H_j(p)$ for some $0\le j\le 4$.
First, we prepare two propositions about 
the structure of the real locus of 1-nodal or 1-cuspidal rational curves.
Both can be readily shown by taking the normalization.

\begin{proposition}\label{prop:nodal}
If $C$ is a rational curve having one ordinary node $p$ as its
only singularity and $C$ is equipped with a real structure $\sigma$,
then the real locus $C^{\sigma}$
consists of one of the following
possibilities:
\begin{itemize}
\item the single point $p$,
\item the union of the single point $p$ and a circle
which does not pass through $p$,
\item the `figure 8' which has a node at $p$.
\end{itemize}
\end{proposition}

\begin{proposition}\label{prop:cuspidal}
If $C$ is a rational curve with a single ordinary cusp $p$ as its only singularity and $C$ is equipped with a real structure $\sigma$, then the real locus $C^{\sigma}$ is a real cuspidal curve homotopic to $\mathbb{S}^1$, with $p$ as its real cusp.
\end{proposition}

To investigate how the structure of the real locus of $S \cap H$ changes as $H$ moves within $(T_pS)^{*\sigma}$, we will use the following proposition.

\begin{proposition}\label{prop:pencil}
For any point $p\in S\us_1\minus\{p_1,p_2\}$,
we have the decomposition
$$
S^{\sigma}\minus \{p\} = \bigsqcup_{H\in (T_pS)^{*\sigma}} \big(S^{\sigma}\minus \{p\}\big)\cap H. 
$$
\end{proposition}

\proof
By the results of the previous subsection, a generic member of the pencil $(T_pS)^{*\sigma}|_S$ is a (rational) curve that has an ordinary node at $p$ as its only singularity.
So if $C$ and $C'$ are distinct generic members of the pencil,
we have $(C,C')_p = 4$ for the local intersection number at $p$.
Since $(C,C')_S = \deg S = 4$, 
this means that the pencil $(T_pS)^*|_S$
has the point $p$ as the unique base point.
Hence, from a basic property of a pencil,
we have
$$
S\minus \{p\} = \bigsqcup_{H\in (T_pS)^*}
\big(S\minus \{p\}\big)\cap H. 
$$
Taking the real locus, we obtain the assertion.
\proofend

\medskip

Our next goal is to determine the adjacent relations for the hyperplanes $H_i(p)$, $0\le i\le 4$.
We already know that these real hyperplanes are mutually distinct
as long as the point $p$ satisfies $\pi_3(p)\not\in B_3^{\sigma}$, and if $\pi_3(p)\in B_3^{\sigma}$, 
then  $H_1(p)=H_3(p)$ as the only coincidence.
By Lemma \ref{lemma:nonmtl3}, if $H\in (T_pS)^{*\sigma}\minus\{H_j(p)\set 0\le j\le 4\}$, then  
$S\cap H$ is a 1-nodal rational curve equipped with a real structure.
So its real locus
is one of the three possibilities in Proposition \ref{prop:nodal}.
The structure of this locus is constant
on each connected component (an open arc) of
the complement $(T_pS)^{*\sigma}\minus\{H_j(p)\set 0\le j\le 4\}$.
We first show:

\begin{lemma}\label{lemma:Hess}
For any $p\in S\us\minus\{p_1,p_2\}$,
the three hyperplanes $H_j(p)$ ($j=2,3,4$) belong to the same  component
of the two arcs $(T_pS)^{*\sigma}\minus\{H_0(p),H_1(p)\}$, 
except for $H_3(p)$ when it is equal to $H_1(p)$.
\end{lemma}

Of course, in the last exceptional situation, $H_3(p)$ is just placed at
the end $H_1(p)$ of the arc whose ends are $H_0(p)$ and $H_1(p)$. (So it is not very exceptional.)

\medskip
\proof
Let $(x,y,z,w)$, $F(x,y),\, G(x,y)$ and ${\bf H}={\bf H}(\lmd,\mu)$ be as in the proof of Lemma \ref{lemma:nonmtl2}.
As in the proof of the same lemma, ${\bf H}(\lmd,\mu)$, regarded as 
a quadratic function on $\CP_1$ has no double root, and  the two roots
correspond to the tangential hyperplane $H_0(p)$ and 
the cuspidal hyperplane $H_1(p)$.
By a linear coordinates change for $(z,w)$, we may suppose that 
$H_0(p)=\{z=0\}$ and $H_1(p) = \{w = 0\}$.
Then the quadratic equation ${\bf H}(\lmd,\mu)=0$ has solutions
precisely at $(\lmd,\mu) = (1,0)$ and $(0,1)$.
Since both solutions are simple, we have either
`${\bf H}(\lmd,\mu)>0$ on the interval $\lmd/\mu>0$ and
${\bf H}(\lmd,\mu)<0$ on the another interval $\lmd/\mu<0$',
or 
`${\bf H}(\lmd,\mu)<0$ on  $\lmd/\mu>0$ and
${\bf H}(\lmd,\mu)>0$ on  $\lmd/\mu<0$'.
By replacing $\lmd$ by $-\lmd$ if necessary, we may suppose that 
the former is the case.

Let $H$ be an element of the real pencil $(T_pS)^{*\sigma}$
defined by the equation \eqref{tpx2} which satisfies $\lmd\mu\neq 0$.
Then since $\bf{H}$ is equal to the discriminant of a quadratic equation of real coefficients,
each of the two branches of the node of $S\cap H$ at $p$
is real if ${\bf H}(\lmd, \mu)>0$, while they are exchanged by 
the real structure $\sigma$ if ${\bf H}(\lmd, \mu)<0$.
So from the above choice of the sign of ${\bf H}$,
the two branches of $S\cap H$ at $p$ are respectively real on the interval $\lmd/\mu>0$
and are exchanged by $\sigma$ on another interval $\lmd/\mu<0$.
On the other hand, by Lemma \ref{lemma:nonmtl1},
the section $S\cap H_j(p)$ ($j=2,3,4$) has an ordinary node at $p$ and the two branches at $p$ are exchanged by $\sigma$, except for $H_3(p)$ when it is equal to $H_1(p)$.
This means that these three hyperplanes belong to 
the interval $\lmd/\mu<0$ under the present choice, apart from the exception.
\proofend

\begin{figure}
\includegraphics{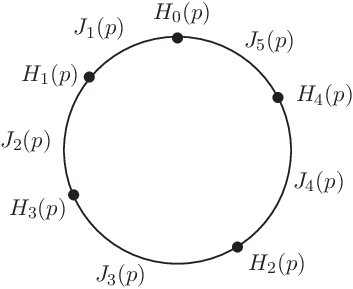}
\caption{
The configuration in the circle $(T_pS)^{*\sigma}$.
}
\label{fig:circle2}
\end{figure}

\medskip
Still let $p$ be a point of $S^{\sigma}_1\minus\{p_1,p_2\}$. 
For a while, we assume $\pi_3(p)\not\in B_3^{\sigma}$,
so that $H_3(p)\neq H_1(p)$.
Then we have five distinct points $H_j(p)$ ($0\le j\le 4$)
on the circle $(T_pS)^{*\sigma}$.
These are the hyperplanes for which their intersections with $S$ are not minitwistor lines (Lemma \ref{lemma:nonmtl3}).
By Lemma \ref{lemma:Hess}, among these five points, the two points $H_0(p)$ and $H_1(p)$ are adjacent
in the circle $(T_pS)^{*\sigma}$.
We define $J_1(p)$ to be the open arc in the circle $(T_pS)^{*\sigma}$ whose ends are $H_0(p)$ and $H_1(p)$.
(So $H_j(p)\not\in J_1(p)$ for any $j=2,3,4$.)
The complement $(T_pS)^{*\sigma}\minus 
\big(J_1(p)\cup\{H_j(p)\set 0\le j\le 4\}\big)$ consists of 
four open arcs.
We name $J_2(p), J_3(p),J_4(p)$ and $J_5(p)$ for these arcs,
where we put the indices in such a way that $J_i(p)$ is adjacent to $J_{i+1}(p)$ where
$J_6(p)$ means $J_1(p)$, and that one of the two ends of $J_2(p)$ is 
the point $H_1(p)$.
(See Figure \ref{fig:circle2}, where the ends of all the arcs are concretely pinned down. They will be determined in the sequel.)

\begin{figure}
\includegraphics{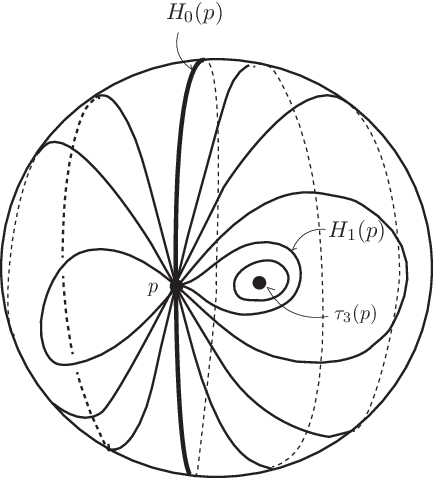}
\caption{
the `foliation' on $S\us_1$ by real hyperplanes in $J_1(p)\cup J_2(p)$
}
\label{fig:Ssigma1}
\end{figure}

\begin{lemma}\label{lemma:5int}
{\em (See Figure \ref{fig:Ssigma1}.)}
Let the situation and notation be as above.
Then the structure of the real locus of a hyperplane section $S\cap H$ is as follows. 
\begin{itemize}
\item[(i)] 
If $H=H_0(p)$, then it is a (double) circle in $S\us_1$.
\item[(ii)] 
If $H\in J_1(p)$, then it is a `figure 8' contained in $S_1^{\sigma}$, which has the point $p$ as the real node of the `figure 8'.
\item[(iii)]
If $H=H_1(p)$, then it is a real cuspidal curve contained in $S_1^{\sigma}$, which has 
the point $p$ as the real cusp.
\item[(iv)]
If $H\in J_2(p)$, then it consists of the point $p$ and
a circle in $S^{\sigma}_1$ which does not pass through $p$.
\end{itemize}
\end{lemma}

\proof
(i) is obvious.
In the following, we use the notations and setting in the proof of Lemma \ref{lemma:Hess}.
Since ${\bf H}(\lmd,\mu)>0$ on the arc $J_1(p)$, when $H\in J_1(p)$,
each of the two branches at the node $p$ of $S\cap H$ is real.
Hence, from Proposition \ref{prop:nodal}, the real locus of
$S\cap H$ has to be `figure 8' having $p$ as the real node. As $p\in S^{\sigma}_1$,
it has to be contained in $S^{\sigma}_1$. So we obtain (ii).
The section $S\cap H_1(p)$ is a cuspidal rational curve that has
the point $p$ as its only singularity, and hence
Proposition \ref{prop:cuspidal} immediately implies the assertion (iii).

When $H\in J_2(p)$, since ${\bf H}(\lmd,\mu)<0$ on $J_2(p)$, 
the two branches of $S\cap H$ at $p$ are exchanged by $\sigma$, 
and hence the node $p$
is an isolated point of $(S\cap H)^{\sigma}$.
To show that $S\cap H$ has a real
(i.e.\,$\sigma$-fixed) circle,
we make use of the real cuspidal curve of $S\cap H_1(p)$.
The intersection of $S$ and $H_1(p)$ is transversal except at the cusp $p$,
and this is the case also for the intersection of the real surface $S^{\sigma}$
and the real hyperplane $H^{\sigma} = \RP_3$.
Because these are moved in a real 4-dimensional space $\RP_4$, 
the last transversality is preserved
under a small displacement of the hyperplane $H_1(p)$ preserving the real structure.
This means that the intersection $S\cap H$ has a real curve
if $H\in J_2(p)$ is sufficiently closed to $H_1(p)$.
Moreover, since the structure of the real locus of $S\cap H$ cannot change when $H$ belongs to the same arc among $J_i(p)$,
the existence of the real circle is valid for any $H\in J_2(p)$.
Then applying Proposition \ref{prop:nodal}, we obtain that 
the real locus of $S\cap H$ consists of the point $p$ 
and a circle that does not go through $p$.
Furthermore, since the real cuspidal curve of $S\cap H_1(p)$ is
contained in $S^{\sigma}_1$, so is the real circle of $S\cap H$
when $H\in J_2(p)$.
Thus we obtain the assertion (iv).
\proofend

\medskip
As above, one of the two ends of the arc $J_2(p)$ is the point $H_1(p)$.

\begin{lemma}\label{lemma:k=3}
Still assuming $H_3(p)\neq H_1(p)$ (namely, $\pi_3(p)\not\in B_3^{\sigma}$),
the other end of the arc $J_2(p)$ is the point $H_3(p)$.
\end{lemma}

\proof
Write $H_i(p)$ ($i=2,3$ or $4$) for the other end of $J_2(p)$.
By Lemma \ref{lemma:nonmtl1}, we have $S^{\sigma} \cap H_i(p) = \{p,\tau_i(p)\}$. 
By Proposition \ref{prop:pencil}, the point $\tau_i(p)$ has to be a shrinking limit of the 
real circle of $S\cap H$ as a point $H\in J_2(p)$ approaches $H_i(p)$.
Since the last circle is contained in the sphere $S^{\sigma}_1$
from Lemma \ref{lemma:5int} (iv),
this means $\tau_i(p)\in S^{\sigma}_1$ by continuity.
Thus, both $p$ and $\tau_i(p)$ belong to the component $S^{\sigma}_1$.
From Proposition \ref{p:symm}, 
the two involutions $\tau_2$ and $\tau_4$ exchange the two spheres
$S^{\sigma}_1$ and $S^{\sigma}_2$, while
$\tau_3$ preserves each of these spheres.
These mean $i=3$.
\proofend

\medskip
Thus, we have obtained the structure of the locus $S^{\sigma}\cap H$
for any $H\in \ol{J_1(p)}\cup\ol {J_2(p)}$,
where $\ol{J_1(p)}$ and $\ol {J_2(p)}$ mean the closures of $J_1(p)$ and $J_2(p)$
in the real circle $(T_pS)^{*\sigma}$ respectively.
When $H=H_0(p)$, the initial situation, the locus $S^{\sigma}_1\cap H$ is a circle, and this circle is 
`double' since $H_0(p)$ is tangent to $S$ in order two along a conic, and the circle $S_1^{\sigma}\cap H$ is exactly the real locus of this conic.
The way how this `double' circle deforms into 
the real cuspidal curve $S^{\sigma}_1\cap H_1(p)$
through the status of `figure 8', while fixing the point $p$  without any overlapping except $p$, is unique at least topologically,
and this transition is illustrated in Figure \ref{fig:Ssigma1}.
Similarly, the transition from the real cuspidal curve $S^{\sigma}_1\cap H_1(p)$ into the two points $S^{\sigma}_1\cap H_3(p)=\{p,\tau_3(p)\}$ through the intermediate status 
$S^{\sigma}_1\cap H\simeq\{p\}\sqcup \mathbb S^1$ ($H\in J_2(p)$) is also unique, and it is as illustrated in Figure \ref{fig:Ssigma1}.
It follows from these that we have
\begin{align}\label{union1}
\bigcup_{H\in \ol {J_1(p)}\,\cup\, \ol{ J_2(p)} } (S\cap H)^{\sigma} = S^{\sigma}_1.
\end{align}

We next discuss the real locus of the nodal curve $S\cap H$ when
a hyperplane $H\in (T_pS)^{*\sigma}$ belongs to the remaining three arcs.
\begin{lemma}\label{lemma:5int2}
As a continuation of Lemmas \ref{lemma:5int} and \ref{lemma:k=3},
\begin{itemize}
\item[(i)] if $H\in J_3(p)$, then the real locus of 
$S\cap H$ consists of the node $p$,
\item[(ii)] if $H\in J_4(p)$, then the real locus of
$S\cap H$ consists of the node $p$ and a circle contained in $S^{\sigma}_2$ (so the circle does not go through $p$),
\item[(iii)] if $H\in J_5(p)$, then the real locus of
$S\cap H$ consists of the node $p$.
\end{itemize}
\end{lemma}

\proof
By Lemma \ref{lemma:k=3}, we have
$\ol {J_2(p)}\cap \ol {J_3(p)} = \{H_3(p)\}$,
and by Lemma \ref{lemma:nonmtl1} we have $S^{\sigma}\cap H_3(p)=\{p,\tau_3(p)\}$.
This is a subset of the sphere $S^{\sigma}_1$,
so $S^{\sigma}_2\cap H_3(p) =\emptyset$.
From the compactness of another sphere $S^{\sigma}_2$,
this means that, if $H\in J_3(p)$ is sufficiently close to the end $H_3(p)$, then $S^{\sigma}_2 \cap H = \emptyset$.
Because the structure of the real locus of $S\cap H$ is unchanged
when $H$ belongs to the same arc, this implies that
we have $S^{\sigma}_2 \cap H = \emptyset$ for any $H\in J_3(p)$.
Since the whole of the sphere $S^{\sigma}_1$ was already swept out
by $H$ belonging to  $\ol {J_1(p)}\cup\ol {J_2(p)}$
as in \eqref{union1}, by Proposition \ref{prop:pencil},
$S^{\sigma}_1\cap H = \{p\}$ for any $H\in (T_pS)^*$ which does not belong to $\ol {J_1(p)}\cup \ol {J_2(p)}$.
So we have $S^{\sigma}\cap H = (S^{\sigma}_1\cup S^{\sigma}_2)\cap H
= (S^{\sigma}_1\cap H)\cup (S^{\sigma}_2\cap H) = \{p\}$
if $H\in J_3(p)$.
Thus we obtain the assertion (i) in the lemma.

Next, we show (iii).
From Lemma \ref{lemma:5int}, if $H=H_0(p)$, then $S\us_2\cap H = \emptyset$. Recalling that $H_0(p)$ is an end of $J_5(p)$, again from the compactness of $S\us_2$, this means that $S\us_2\cap H = \emptyset$ if $H\in J_5(p)$.
Using \eqref{union1}, this means (iii).

\begin{figure}
\includegraphics{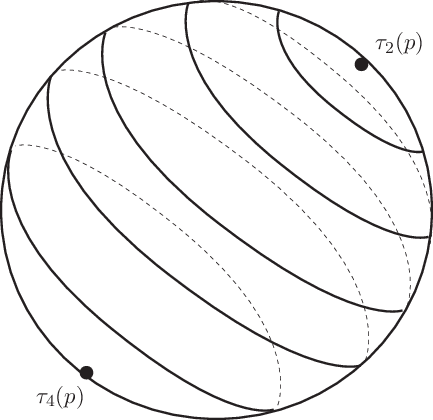}
\caption{
the `foliation' on $S\us_2$ by real hyperplanes in $J_4(p)$
}
\label{fig:Ssigma2}
\end{figure}

For the remaining assertion (ii), 
since any $H\in  {J_5(p)}$ does not intersect $S\us_2$ from what we have obtained, it follows from Proposition \ref{prop:pencil} that  
\begin{align}\label{union2}
\bigsqcup_{H\in \ol {J_4(p)} } S^{\sigma}_2\cap H = S^{\sigma}_2.
\end{align}
From Proposition \ref{prop:nodal}, this means that 
$S^{\sigma}_2\cap H$ is either always a circle or always the figure 8 when $H\in J_4(p)$.
Let $H_j(p)$ and $H_k(p)$ the endpoints of $J_4(p)$, so that $\{j,k\} = \{2,4\}$.
Later we will show $(j,k) = (2,4)$ (Remark \ref{rmk:H234}), but at this stage we do not need it.
By Lemma \ref{lemma:nonmtl1},  $S^{\sigma}\cap H_j(p) = \{p,\tau_j(p)\}$ and $S^{\sigma}\cap H_k(p) = \{p,\tau_k(p)\}$, and from Proposition \ref{p:symm}, the points $\tau_j(p)$ and $\tau_k(p)$ belong to  $S^{\sigma}_2$. 
From \eqref{union2}, this means that 
\begin{align}\label{union4}
\bigsqcup_{H\in{J_4(p)} } S^{\sigma}_2\cap H = S^{\sigma}_2
\minus \{\tau_j(p),\tau_k(p)\}.
\end{align}
Obviously this is impossible if the curve $S^{\sigma}_2\cap H$ is a `figure 8'.
Hence $S^{\sigma}_2\cap H$ has to be a circle and 
it shrinks to $\tau_2(p)$ or $\tau_4(p)$ as $H$ approaches the endpoints of $J_4(p)$. (See Figure \ref{fig:Ssigma2}.)
Thus, we have obtained assertion (ii) of the lemma.
\proofend

\medskip
Next, we discuss the case where the real point $p\in S^{\sigma}_1\minus\{p_1,p_2\}$ satisfies $\pi_3(p)\in B_3^{\sigma}$, so that 
$H_3(p) = H_1(p)$.
We note that while $B_3^{\sigma}$ and hence $\pi_3\inv(B_3^{\sigma})$ also consist of two mutually disjoint circles,
the intersection $S^{\sigma}_1\cap \pi_3\inv(B_3^{\sigma})$ consists
of a single circle.
The four hyperplanes $H_0(p), H_1(p), H_3(p)$ and $H_4(p)$
are all ones in $(T_pS)^{*\sigma}$ whose intersections with $S$
are not minitwistor lines.
We name the four open arcs of $(T_pS)^{*\sigma}\minus 
\{H_i(p)\set i=0,1,3,4\}$ $J_1(p),J_3(p),J_4(p)$ and $J_5(p)$
in such a way that these are arranged in this order,
the two ends of $J_1(p)$ are $H_0(p)$ and $H_1(p) \,(=H_3(p))$,
and that one of the ends of $J_3(p)$ is $H_1(p)$.
Here, we skipped $J_2(p)$ because we can regard the arc
$J_2(p)$ appearing in the case $\pi_3(p)\not\in B_3^{\sigma}$
to shrink to one point $H_1(p)=H_3(p)$ as a point $p\in S^{\sigma}_1$ 
approaches a point of $S^{\sigma}_1\cap\pi_3\inv(B_3^{\sigma})$.

\begin{lemma}\label{lemma:5int3}
Under the above situation and notations,  
\begin{itemize}
\item[(i)] 
if $H\in J_1(p)$, then the real locus of $S\cap H$ is a `figure 8' contained in $S^{\sigma}_1$, which has the point $p$ as the real node,
\item[(ii)] if $H\in J_3(p)\cup J_5(p)$, then the real locus of $S\cap H$ consists
of the node $p$,
\item[(iii)] If $H\in J_4(p)$, then the real locus of $S\cap H$ consists of the node $p$ and a circle contained in $S^{\sigma}_2$.
\end{itemize}
\end{lemma}

\proof
The proof for  (ii) in Lemma \ref{lemma:5int} works without any change for the present assertion (i).
As $H_1(p) = H_3(p)$, and $\pi_3(p)\in B_3$
we have $S^{\sigma}\cap H_1(p)=\{p=\tau_3(p)\}$, and by continuity, the `figure 8' {\em directly} shrinks to the point $\tau_3(p)$ without passing a real cuspidal curve as $H\in J_1(p)$ approaches the point $H_3(p)$. Therefore, slightly different from \eqref{union1}, this time we have
\begin{align}\label{union3}
\bigcup_{H\in \ol {J_1(p)}} (S\cap H)^{\sigma} = S^{\sigma}_1.
\end{align}
Hence, from Proposition \ref{prop:pencil}, if 
$H\in (T_pS)^{*\sigma}\minus \ol{J_1(p)}$, then $H$ cannot intersect
$S^{\sigma}_1\minus\{p\}$.
Further, by the same argument as the case (i) in Lemma \ref{lemma:5int2}, we obtain that $S^{\sigma}_2\cap H = \emptyset$ if $H\in J_3(p)$. Together with \eqref{union3}, this means  $S^{\sigma}\cap H =\{p\}$ if $H\in J_3(p)$.
Hence we obtain (ii) for the case $H\in J_3(p)$.
The assertions in the case $H\in J_5(p)$ and (iii) can be shown in the same way as
(iii) and (ii) in Lemma \ref{lemma:5int2} respectively.
\proofend

\medskip We remark that, from the proof, the equality \eqref{union2}
still holds in the case $\pi_3(p)\in B_3^{\sigma}$ too.
Namely, the equality \eqref{union2} is valid for any point $p\in S^{\sigma}_1\minus\{p_1,p_2\}$.

\subsection{Selecting a component from the space of real minitwistor lines}\label{ss:rmtl}
As we have assumed, let $S\subset\CP_4$ be the Segre surface defined by 
\eqref{Segre1}, and $S^*\subset\CP_4^*$ the projectively dual variety of $S$.
Let $I(S)\subset\CP_4\times\CP_4^*$ be the incidence variety of $S$,
which is by definition the closure of 
the locus 
$$
\big\{(p,H)\in S_{\reg}\times\CP_4^*\set 
T_pS\subset H\big\}
$$
taken in $\CP_4\times\CP_4^*$.
There is a double fibration
\begin{align}\label{diagram:double1}
S \stackrel{\pr_1}\longleftarrow I(S) \stackrel{\pr_2}\longrightarrow S^*
\end{align}
where $\pr_1$ and $\pr_2$ denote the projections to the
factors respectively.
The real structure $\sigma$ on $S$ (defined by \eqref{rs1})
naturally induces those on $S^*$ and $I(S)$.
We use the same symbol $\sigma$ for these.
Using the component $S^{\sigma}_1$ of $S^{\sigma}$, we define
\begin{align}\label{ISsigma0}
I(S)^{\sigma}_1:=\pr_1\inv(S^{\sigma}_1) \cap I(S)^{\sigma}.
\end{align}
Since $S^{\sigma}_1\subset S_{\reg}$, this is a real submanifold of $I(S)$,
and the projection $\pr_1$ gives 
the structure of an $\mathbb S^1$-bundle
\begin{align}\label{ISsigma}
I(S)^{\sigma}_1 \lras S^{\sigma}_1 
\end{align}
whose fiber over a point $p\in S^{\sigma}_1$
may be identified with the real pencil $(T_pS)^{*\sigma}\simeq \mathbb S^1$.

Next, we define
$$
S^{*\sigma}_1 := \pr_2\big(I(S)^{\sigma}_1\big).
$$
This is the space of real hyperplanes in $\CP_4$ which are tangent to $S$ at a point of $S\us_1$.
Then by restricting \eqref{diagram:double1} to $I(S)^{\sigma}_1$, we obtain the double fibration
\begin{align}\label{diagram:double2}
S\us_1 \stackrel{\pr_1}\longleftarrow  I(S)\us_1  
\stackrel{\pr_2}\longrightarrow S^{*\sigma}_1
\end{align}

We note that the $\CC^*$-action on $S$ naturally induces that on $I(S)$ and $S^*$ and the diagram \eqref{diagram:double1} is $\CC^*$-equivariant.
Similarly, the diagram \eqref{diagram:double2} is $\mathbb S^1$-equivariant.
The $\mathbb S^1$-action on $S\us_1$ is the rotation that fixes the two points $p_1$ and $p_2$, and points on the fibers over these points are fixed by the action.

To each index $j=0,2,3,4$,
the mapping $p\longmapsto (p, H_j(p))$ for $p\in S^{\sigma}_1$ gives a 
smooth section of \eqref{ISsigma} which is globally defined on $S^{\sigma}_1$.
We denote $\ms H_j$ for the images of these four sections.
Because $S^{\sigma}_1\simeq \mathbb S^2$ is simply connected,
these sections  
divide the $\mathbb S^1$-bundle $I(S)^{\sigma}_1$ into four connected components, and each component is diffeomorphic to 
the trivial interval bundle over $S^{\sigma}_1$.

We are interested in the connected component bounded by
the two sections $\ms H_2$ and $\ms H_4$.
In the sequel, we always denote this component by the letter $W$,
and its compactification $W\cup \ptl W$ where $\ptl W = \ms H_2\sqcup\ms H_4$ by $\ol W$. We call $\ms H_2$ and $\ms H_4$ the {\em future infinity} and the {\em past infinity} of $W$. 
Further, we denote the restriction of the projection \eqref{ISsigma} to $W$ and $\ol W$ by $\varpi$ and $\ol\varpi$ respectively, so that we have the projections
\begin{align}\label{W}
\varpi:W \lras S^{\sigma}_1
\qandq
\ol\varpi:\ol W\lras S\us_1.
\end{align}
These are fiber bundle map whose fiber over a point $p\in S\us_1$ is the open arc $J_4(p)$ and its closure $\ol {J_4(p)}$, and as above, 
these are trivial bundles over $S\us_1\simeq\mathbb S^2$.
So letting $I$ be an open interval in $\RR$,
$W$ is diffeomorphic to $\mathbb S^2\times I$ and $\ms H_2\sqcup\ms H_4$ are its boundary.

\begin{proposition}\label{prop:EW3}
The map $\pr_2|_W$ from $W\simeq\mathbb S^2\times\RR$
to the image $\pr_2(W)\subset S^*$ is a real analytic diffeomorphism.
\end{proposition}

\proof
The real analyticity  of $\pr_2|_W:W\to\pr_2(W)$ is immediate since $\pr_2:I(S)\to S^*$ is a holomorphic map 
and $W$ is a real slice of (a part of) the smooth locus of $I(S)$.
To show diffeomorphicity of $\pr_2|_W$, from Proposition \ref{prop:2sing},
the fiber $\pr_2\inv(H)$ of a point $H\in S^*$ has at least two points whose images to $S$ by $\pr_1$ are smooth points of $S$
only if $H$ belongs to 
\begin{align}\label{noninj}
Q_2^*\cup Q_3^*\cup Q_4^*\cup \Lmd^*\cup\big\{ f\inv(\ol{\lmd_i\lmd_j})\set 1\le i<j\le 4\big\},
\end{align}
where $Q_j^*$ ($j=2,3,4)$ denote the images of
the dual quadrics $Q_j^*\subset\CP_3^*$ by the dual maps of the projections $\pi_j:\CP_4\lras \CP_3$, and $\Lmd^*$ denotes the image of
the dual conic $\Lmd^*\subset\CP_2^*$ by the dual map of the projection $f:\CP_4\lras \CP_2$.
From the definition of the space $W$, it is easy to see that 
the image $\pr_2(W)$ is disjoint from the set \eqref{noninj}.
Therefore, if $A$ denotes the inverse image of the set \eqref{noninj} to $I(S)$ under $\pr_2$, then $W\cap A = \emptyset$.
Further, $A$ is a closed subset of $I(S)$ since the set \eqref{noninj} is a closed subset of $S^*$.
We do not know if $\pr_2$ is injective on the whole of $I(S)\minus A$,
but it is injective on $I(S)\minus 
(A\cup \pr_1\inv\{e_0,e_1\})$ from Proposition \ref{prop:2sing}.
Further, from the definition of $W$, we have $ \pr_1\inv\{e_0,e_1\}\cap W = \emptyset$. 
Therefore, $W\subset I(S)\minus 
(A\cup \pr_1\inv\{e_0,e_1\})$.
Further, $I(S)\minus 
(A\cup \pr_1\inv\{e_0,e_1\})$ is open in $I(S)$.
Therefore, $I(S)\minus 
(A\cup \pr_1\inv\{e_0,e_1\})$ is an open neighborhood of $W$ in $I(S)$ on which $\pr_2$ is injective.
Since any injective holomorphic map between complex manifolds of the same dimension is holomorphic onto its image
\cite[p.\,19,\,Proposition]{GrHr}, this means that in a neighborhood of $\pr_2$ is biholomorphic onto its image.
Therefore, its restriction to $W$ is diffeomorphic onto its image.
\proofend

\medskip
In the following, we identify $W$ with its image $\pr_2(W)\subset S^*$. 

We define a dense open subset of $W$ by
\begin{align}\label{W0}
W_0&:= \varpi\inv(S\us_1\minus\{p_1,p_2\})
=W\minus \big(\varpi\inv(p_1)\cup \varpi\inv(p_2)\big).
\end{align}
Based on the result from the previous section, this is the space of real minitwistor lines that have a node at a point on $S^\sigma_1$ and intersect $S^\sigma_2$ in a circle.
From this and a result in \cite{HN11}, we readily obtain 

\begin{proposition}\label{prop:W0}
The 3-manifold $W_0$ admits a Lorenzian Einstein-Weyl structure.
\end{proposition}

\proof
In \cite[Theorem 4.3]{HN11}, we showed that the space of real minitwistor lines has a positive definite EW structure when the minitwistor lines have only isolated real points. This proof also applies when the real minitwistor lines have real circles; in such cases, the induced EW structure becomes indefinite.
\proofend

\medskip
The following characterization of 
the manifold $W$ is evident from the construction 
(and Lemmas \ref{lemma:5int}--\ref{lemma:5int3}).

\begin{proposition}\label{prop:charW}
The 3-manifold $W$ in \eqref{W} is exactly the space of 
all real regular minitwistor lines
and all real irregular minitwistor lines in $S$,
such that they have a node at a point of $S^{\sigma}_1$
and such that they intersect $S^{\sigma}_2$ in a circle.
The regular ones are parameterized by the open subset $W_0$ of $W$ and 
the irregular ones are parameterized by the two fibers $\varpi\inv(p_1)$ and $\varpi\inv(p_2)$.
\end{proposition}

\begin{remark}\label{rmk:2pts}
{\em
According to Lemmas \ref{lemma:5int}--\ref{lemma:5int3}, the statements of the proposition remain valid if we replace the condition ``such that they intersect $S^{\sigma}_2$ in a circle'' with ``such that they intersect $S^{\sigma}_2$ at {\em at least two points}''. This stronger version will be used in a crucial part of the proof for Zoll property of our EW structures (Theorem \ref{thm:Zoll}), specifically in Proposition \ref{prop:slg2}.}
\end{remark}

\subsection{The irregular minitwistor lines}\label{ss:imtl}
Our next goal is to show that the EW structure on $W_0$ obtained in the previous section extends smoothly to the two intervals $\varpi\inv(p_i)$ ($i=1,2$),
the $\mathbb S^1$-fixed locus of $W\simeq\mathbb S^2\times I$.
In this subsection, for this purpose, we discuss the hyperplane sections by $H$ belonging to the two intervals $\varpi\inv(p_i)$ ($i=1,2$).
If $H$ is such a real hyperplane, then since the tangent space $T_{p_i}S$ is spanned by the two lines $l_i$ and $\ol l_i$,
the hyperplane section $S\cap H$ is of the form $l_i + \ol l_i + D$ for some conic $D$. Since $e_0\in l_i$ and $e_1\in \ol l_i$,  $H$ is of the form $f\inv(l)$
for some real line $l\subset\CP_2$ and it passes through the discriminant point $\lmd_i$.
Letting $\lmd\in\Lmd$ be the residual point of $\Lmd\cap l$,
this means that the conic $D$ is of the form $f\inv(\lmd)$, and $\lmd\in\Lmd\us$.
We denote $D_{\lmd}$ for this fiber conic.
Thus, if $H\in \varpi\inv(p_i)$ ($i=1,2$), then 
\begin{align}\label{imtl1}
S\cap H = l_i + \ol l_i + D_{\lmd}
\end{align}
for some $\lmd\in \Lmd\us$.
Therefore, there is a natural inclusion $\varpi\inv(p_1)\cup \varpi\inv(p_2)\subset \Lmd\us$. 
Under this understanding, the intervals $\varpi\inv(p_1)$ and $\varpi\inv(p_2)$ are indentied with the arc $I_2$ (over which the sphere $S\us_2$ lie). More precisely, we have:

\begin{proposition}\label{prop:2intervals}
Under this understanding $\varpi\inv(p_i)\subset \Lmd^{\sigma}$ for $i=1,2$, both of the intervals are identified with the arc $I_2 = (\lmd_3,\lmd_4)\subset \Lmd^{\sigma}$ defined in Definition \ref{def:Ssigma}.
Under this identification, the ends $H_2(p_1)$ and $H_4(p_1)$ of $\varpi\inv(p_1)$
are respectively correspond to the ends $\lmd_3$ and $\lmd_4$ of $I_2$,
while the ends $H_2(p_2)$ and $H_4(p_2)$ of $\varpi\inv(p_2)$ are respectively correspond to the ends $\lmd_4$ and $\lmd_3$ of $I_2$.
\end{proposition}

\proof
Using the involution $\tau_2$ of $S$ in Section \ref{s:S}, 
we readily obtain $\pi_2 (p_1) = \pi_2 (p_3)$ for the double covering map $\pi_2:S\lras Q_2$. This means that 
\begin{align*}
S\cap H_2(p_1) &= \big(l_1 + \ol l_1\big) + \big(l_3 + \ol l_3\big).
\end{align*}
As this equals $S\cap f\inv (\ol{\lmd_1\lmd_3})$,
$H_2(p_1) = f\inv (\ol{\lmd_1\lmd_3})$.
Hence, the end $H_2(p_1)$ of $J_4(p_1)$ corresponds to the discriminant point $\lmd_3$.
Similarly, we obtain 
$H_2(p_2) = f\inv (\ol{\lmd_2\lmd_4}),$ 
$H_3(p_1) = H_3(p_2) = f\inv (\ol{\lmd_1\lmd_2})$,
$H_4(p_1) = f\inv (\ol{\lmd_1\lmd_4})$ and 
$H_4(p_2) = f\inv (\ol{\lmd_2\lmd_3})$.
From these, 
we obtain the assertions of the proposition.
\proofend

\begin{remark}\label{rmk:H234}{\em
From the explicit forms of the hyperplanes $H_j(p_i)$
for $i=1,2$ and $j=2,3,4$ obtained in the proof of Proposition \ref{prop:2intervals}, we can conclude $\ol {J_3(p)}\cap\ol {J_4(p)} = H_2(p)$ and $\ol {J_4(p)}\cap\ol {J_5(p)} = H_4(p)$.
}
\end{remark}

Now we give a definition, which will be important in the rest of this paper.

\begin{figure}
\includegraphics{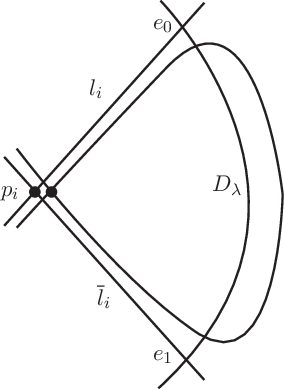}
\caption{
an irregular minitwistor line and a regular minitwistor line
}
\label{fig:fan}
\end{figure}

\begin{definition}\label{d:imtl}
{\em (See Figure \ref{fig:fan}.)
For any point $\lmd\in I_2 = (\lmd_3,\lmd_4)\subset\Lmd\us$, we call the real reducible curve $l_i + \ol l_i + D_{\lmd}$ in \eqref{imtl1} an {\em irregular minitwistor line}. In other words, these are sections by hyperplanes in $\varpi\inv(p_1)\cup\varpi\inv(p_2)$. 
For distinction, we call the hyperplane sections corresponding to any point of $W_0$ {\em regular minitwistor lines}.\proofend
}
\end{definition}

This terminology indicates that, despite their reducibility, irregular minitwistor lines still play a role similar to that of true minitwistor lines, as they determine points on the EW manifold. Note that all irregular minitwistor lines are $\mathbb{C}^*$-invariant.

Next, to discuss hyperplane sections that correspond to points on the future and past infinity $\ms H_2$ and $\ms H_4$, for $j\in\{2,4\}$,
recall that there are double covering maps $\pi_j:S\lras Q_j$ to smooth quadric $Q_j$ that preserve the real structure.
Then using Proposition \ref{prop:3rs}, we readily see that for any point $p\in S\us_1$, we can write
\begin{align}\label{AolA}
S\cap H_j(p) = \pi_j\inv(\CCC) + \pi_j\inv(\ol \CCC),
\end{align}
where $\CCC$ and $\ol\CCC$ are the pair of lines on $Q_j$ that intersect at the point $\pi_j(p)$.
So $\CCC$ and $\ol\CCC$ are $(1,0)$- and $(0,1)$-curves on $Q_j\simeq\qdr$.
The curves $\pi_j\inv(\CCC)$ and $\pi_j\inv(\ol \CCC)$ on $S$ are conics exchanged by $\sigma$ and these intersect transversely at the two points $p\in S\us_1$ and $\tau_j(p)\in S\us_2$.
These conics split into two lines only when 
the point $p\in S\us_1$ is either $p_1$ or $p_2$.
The decomposition \eqref{AolA} will be needed in the next section.

\section{The properties of the EW spaces}\label{s:EW}

\subsection{The modification of the Segre surfaces}
\label{ss:modify}
In this subsection, we introduce a modification to a part of any Segre surface $S$. This modification will be used to demonstrate the extension property of the EW structure on $W_0$ across the two intervals and to analyze the variation of EW structures as the complex structure of $S$ undergoes deformation.

\begin{figure}
\includegraphics{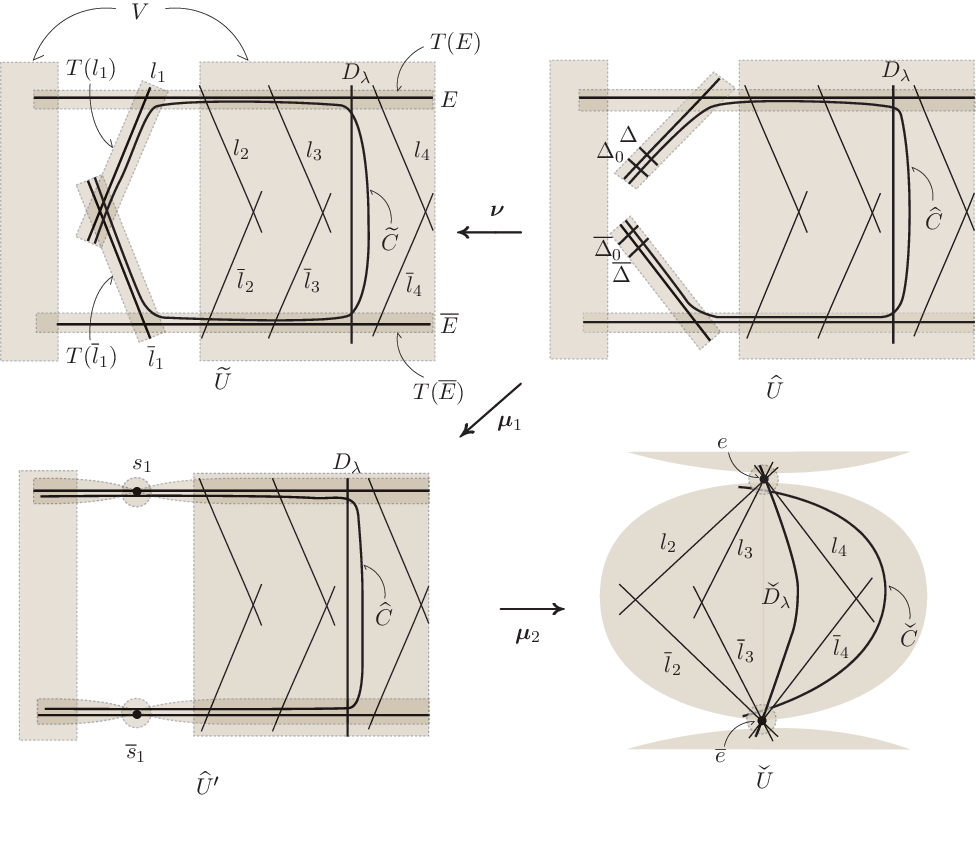}
\caption{
The open surface $U$ and its modification $\check U$
}
\label{fig:modif1}
\end{figure}

In the following, we first choose an open real connected subset $U$ of $S$, which can be thought as a tubular neighborhood of any irregular minitwistor line of the form $l_1 + \ol l_1 + D_{\lmd}$ ($\lmd\in I_2$).
Let $\mu:\tilde S\lras S$ be the minimal resolution of the two nodes $e_0$ and $e_1$ of $S$, and $E$ and $\ol E$ the exceptional curves of $e$ and $\ol e$ respectively. These are $(-2)$-curves on $\tilde S$.
Then $\mu$ eliminates the indeterminacy of the projection $f:S\lras \Lmd$. We write $\tilde f:\tilde S\lras\Lmd$ for the resulting holomorphic map.
Take a real (i.e.\,$\sigma$-invariant) open disk $\mathbb D\subset\Lmd$ 
that satisfies $\lmd_1\in \mathbb D$ and $\lmd_2,\lmd_3,\lmd_4\not\in\ol{\mathbb D}$, where $\ol{\mathbb D}$ denotes the closure of $\mathbb D$.
Define $V := \tilde f\inv\big( \Lmd \minus\ol{\mathbb D}\big)$.
This is a real $\mathbb S^1$-invariant open subset of $\tilde S$ and  satisfies $V\supset \cup_{i=2,3,4}(l_i\cup\ol l_i)$. 
The set $\tilde S \minus \ol V$ is an open neighborhood of the fiber $\tilde f\inv(\lmd_1) = l_1 + \ol l_1$.
Take a Riemannian metric on $\tilde S$ that is invariant under the $\mathbb S^1$-action and the real structure and 
let $T(l_1)$ be a tubular neighborhood of $l_1\subset\tilde S$ with respect to the metric.
Let $T(\ol l_1)$ be the image of $T(l_1)$ under $\sigma$. This is an open tubular neighborhood of $\ol l_1$. 
These two neighborhoods are also $\mathbb S^1$-invariant.
By choosing the tube sufficiently small, we may suppose that 
$T(l_1)\cup T(\ol l_1)
\subset
\tilde S\minus \overline V.$
Finally, we choose a tubular neighborhood $T(E)$ of $E$ using the above Riemannian metric on $\tilde S$ in a way that 
$T(E)\cap T(\ol E)=T(E)\cap T(\ol l_1) = \emptyset$, where
$T(\ol E):=\sigma(T(E))$.
These are also $\mathbb S^1$-invariant.
Using these open subsets, we define subsets $\tilde U\subset\tilde S$ 
and $U\subset S$ by
\begin{align}\label{tSe}
\tilde U:= V\cup T(l_1)\cup T(\ol l_1) \cup 
T(E) \cup T(\ol E) 
\qandq
U:= \mu\big(\tilde U\big).
\end{align}
(See Figure \ref{fig:modif1} for $\tilde U$.)
These are real connected $\mathbb S^1$-invariant open subsets, where
openness of $U$ follows from $\tilde U\supset T(E)\cup T(\ol E)$.

Next, starting from the complex surface $\tilde U$, we define a series of operations.
Firstly, consider the disjoint union
$\hat T(l_1,\ol l_1):= T(l_1)\sqcup T(\ol l_1).$
This is a disconnected complex surface equipped with
an $\mathbb S^1$-action and a real structure that flips the two connected components $T(l_1)$ and $T(\ol l_1)$.
Then we can define 
a non-singular complex surface $\hat U$ by replacing 
$T(l_1)\cup T(\ol l_1)$ in \eqref{tSe} with 
$\hat T(l_1,\ol l_1)$.
$\hat U$ is a connected complex surface which has a real structure and an $\mathbb S^1$-action.
We write 
\begin{align}\label{}
\boldsymbol\nu: \hat U \lras \tilde U
\end{align}
for the natural surjection. This is a holomorphic map preserving the real structure and the $\mathbb S^1$-action.
The inverse image $\boldsymbol\nu\inv(l_1)$ consists of two connected components and one of them is the component in $T(l_1)\subset \hat U$, which is mapped onto $l_1$ by $\boldsymbol\nu$.
We use the same letter $l_1$ for this component.
Another component is a piece of smooth curve in $T(\ol l_1)\subset \hat U$ that is mapped by $\boldsymbol\nu$ onto the portion $l_1\cap T(\ol l_1)$.
(See Figure \ref{fig:modif1}.)
We denote this piece of curve by $\ol\DDD_0$.
This is biholomorphic to a disk in $\CC$.
We then have $\boldsymbol\nu\inv(l_1) = l_1\sqcup \ol \DDD_0$ and $\boldsymbol\nu\inv(\ol l_1) = \ol l_1 \sqcup \DDD_0$. 
These pieces of curves will be significant later.

Since $\boldsymbol\nu$ gives a biholomorphic map from $T(l_1)\subset \hat U$ onto $T(l_1)\subset \tilde U$,
the self-intersection numbers of $l_1$ in $\tilde U$ and $\hat U$ are equal. So $l_1^2 = -1$ on $\hat U$. Similarly, $\ol l_1^2 = -1$ on $\hat U$.
As a second step of the modification, let 
$
\boldsymbol\mu_1: \hat U \lras \hat U'
$
be the blowing-down of these two $(-1)$-curves.
We denote $s_1$ and $\ol s_1$ for the points 
$\boldsymbol\mu_1(l_1)$ and $\boldsymbol\mu_1(\ol l_1)$ respectively and for simplicity we still write $E$ and $\ol E$ for the images $\boldsymbol\mu_1(E)$ and $\boldsymbol\mu_1(\ol E)$ respectively, so that $s_1\in E$ and $\ol s_1\in \ol E$.
Since $l_1\sqcup\ol l_1\subset\hat U$ are real and $\mathbb S^1$-invariant, the surface $\hat U'$ has an $\mathbb S^1$-action and the real structure induced from those on $\hat U$.

Evidently, the self-intersection numbers of $E$ and $\ol E$ change through $\boldsymbol\mu_1$ and $E^2 = \ol E^2 = -1$ on the surface $\hat U'$.
As a final step of the modification, let 
$
\boldsymbol\mu_2: \hat U' \lras \check U
$
be the blowing-down of these two $(-1)$-curves.
Again, this preserves the real structure and the $\mathbb S^1$-action.
We still denote $\sigma$ for the real structure on $\check U$ and put $e:=\boldsymbol\mu_2(E)$
and $\ol e:= \boldsymbol\mu_2(\ol E)$.
These two points are fixed by the $\mathbb S^1$-action.
For any $2\le i\le 4$, we use the same letters $l_i$ and $\ol l_i$ for the transformations into $\check U$ of the six lines in $S$. 
In $\check U$, there do not exist $l_1$ and $\ol l_1$. 
As an effect of $\boldsymbol\mu_2$, $l_i^2 = \ol l_i^2 = 0$ on 
$\check U$.
The blowing-down $\boldsymbol\mu_2$ terminates the modification.
Considering the open subset $U = \mu(\tilde U)$ of $S$
as in \eqref{tSe}, we can regard $\check U$ as a transformation of $U$.

We emphasize that we cannot apply this operation to the entire surface $\tilde S$.
We also note that the whole operation works when we replace the role of $l_1+\ol l_1$ by $l_2+\ol l_2$.

\subsection{The transformations of minitwistor lines}\label{ss:extension}
The subset $U$ in the previous subsection contains all irregular minitwistor lines of the form 
$l_1 + \ol l_1 + D_{\lmd}$, where $\lmd\in I_2$.
On the other hand, the inverse images of regular minitwistor lines in $S$ into $\tilde S$ are of degree two over $\Lmd$,
and therefore, they are in general not contained in the open surface $\tilde U$ as they can across the `blank area' $\tilde S\minus \tilde U$.
But if a regular minitwistor line is sufficiently close to the irregular ones (i.e.\,if the regular ones belong to a sufficiently small neighborhood of the interval $\varpi\inv(p_1)$), then it is contained in $U$.
In this subsection, we show that all these minitwistor lines in $U$ are transformed into {\em smooth} minitwistor lines in $\check U$.

First, take any point $\lmd\in I_2$ and write $\check D_{\lmd}$
for the transformation of the fiber conic $D_{\lmd}\subset U$ (or $l_1+\ol l_1 + D_{\lmd}\subset \tilde U$) into $\check U$.
We easily see that  $e,\ol e\in \check D_{\lmd}$, $\check D_{\lmd}^2 = 2$, and that
all these transformations are linearly equivalent.
In particular, all curves $\check D_{\lmd}\subset \check U$ are smooth minitwistor lines.
All these are $\mathbb S^1$-invariant.

Next, we investigate the transformations of regular minitwistor lines in $U$.
Let $C\subset U$ be such a one and write $\tilde C:=\mu\inv(C)$.
The inverse image $\boldsymbol\nu\inv(\tilde C)\subset \hat U$ consists of three irreducible components.
One of them is mapped onto $\tilde C$ by $\boldsymbol\nu$ and this is a {\em smooth} rational curve.
We denote $\hat C$ for this component. (See Figure \ref{fig:modif1}.) This can be regarded as the normalization of the nodal rational curve $C$ or $\tilde C$.
We denote $\DDD$ and $\ol\DDD$ for the remaining two irreducible components of $\boldsymbol\nu\inv(\tilde C)$.
Similarly to $\DDD_0$ and $\ol \DDD_0$ that arose from $l_1$ and $\ol l_1$, these are pieces of smooth curves.
We make a distinction between $\DDD$ and $\ol\DDD$ by supposing that $\DDD$ intersects $l_1$. (See Figure \ref{fig:modif1}.)
Since the map $\hat C\lras C$ is the normarization of the node of $C$, we have $\hat C^2 = 4-2 = 2$ on $\hat U$.
Since $\hat C\subset\hat U$ intersects none of $l_1,\ol l_1,E$ and $\ol E$, the complex structure of a neighborhood of $\hat C$ in $\hat U$ is not affected by the blowdowns $\boldsymbol\mu_1$ and $\boldsymbol\mu_2$, and if we denote $\check C:=\boldsymbol\mu_2\circ\boldsymbol\mu_1(\hat C)$, then $\check C^2=2$ on $\check U$.
Thus, regular minitwistor lines in $U$ are transformed again into smooth minitwistor lines in $\check U$.
Obviously, these minitwistor lines pass through none of the two $\mathbb S^1$-fixed points $e$ and $\ol e$.
These minitwistor lines are not $\mathbb S^1$-invariant.

Summing up, we have obtained:

\begin{proposition}\label{p:modmtl}
The modification of the open subset $U\subset S$ to $\check U$ given in the previous subsection transforms
regular and irregular minitwistor lines in $U$ into smooth minitwistor lines in $\check U$.
More precisely, irregular minitwistor lines are transformed into the ones passing through the two points $e$ and $\ol e$, while regular ones are transformed into the ones not passing through $e$ nor $\ol e$.
\end{proposition}

From this, we obtain the required extension property of the EW structure:

\begin{theorem}\label{t:WEW}
The EW structure on $W_0$ in Proposition \ref{prop:W0} extends real analytically across the two intervals $\varpi\inv(p_i)$ ($i=1,2$). 
\end{theorem}

\noindent{\em Proof.}
By construction, the minitwistor space $\check U$ contains the transformations of all irregular minitwistor lines of the form $l_1 + \ol l_1 + D_{\lmd}$, where $\lmd\in I_2$, as well as the transformations of all nearby regular minitwistor lines in $U$.
This means that the EW structure on $W_0$ extends real analytically across the interval $\varpi\inv(p_1)$.
The same conclusion can be deduced for another interval 
$\varpi\inv(p_2)$
by just exchanging the role of $l_1+\ol l_1$ and $l_2+\ol l_2$ (including the modification construction).
\proofend

\medskip
Thus, from the Segre surface $S$, we have obtained the real analytic Lorenzian EW structure on the 3-manifold $W=W_0\cup\varpi\inv(p_1)\cup\varpi\inv(p_2)\simeq\mathbb S^2\times I$.

\begin{remark}{\em
Originally, we proved the extension property as stated in Theorem \ref{t:WEW} by employing the deformation theory of singular curves on a surface (which we will briefly discuss in Section \ref{ss:slg}), not using the modification techniques.
 This approach has the advantage that 
 the dual variety 
$S^*$ is known to be smooth even at points corresponding to irregular minitwistor lines, for example.
}\end{remark}

As above, all minitwistor lines on $\check U$ are smooth and the distinction for those coming from regular or irregular ones on $U$ is reflected in whether they pass through $e$ and $\ol e$, or equivalently, the $\mathbb S^1$-invariance.
The next proposition gives another interesting distinction.

\begin{proposition} \label{prop:le1} 
On the open surface $\check U$,
the $\mathbb S^1$-invariant minitwistor lines are mutually linearly equivalent and any such a minitwistor line $\check D_{\lmd}$ satisfies $\dim |\check D_{\lmd}| = 1$.
If a minitwistor line $\check C\subset \check U$ is obtained from a regular minitwistor line $C\subset U$,
then $\dim|\check C| = 0$.
\end{proposition}

\proof
We only prove the case that $\check U$ is obtained from $U$ that contains $l_1+\ol l_1$ since the case $l_2+\ol l_2\subset U$ can be obtained similarly.

Any minitwistor lines on $S$ are linearly equivalent to each other since they are all obtained as hyperplane sections of $S\subset\CP_4$.
Hence, any minitwistor lines on the open subset $U\subset S$ are also linearly equivalent.
As is already remarked, the linear equivalence $D_{\lmd}\sim D_{\lmd'}$ for any $\lmd,\lmd'\in I_2$ follows immediately from this and our modification procedure.

To show $\dim |D_{\lmd}|=1$ on $\check U$, let $C_0 := l_1 + \ol l_1 + D_{\lmd}$ 
($\lmd\in I_2$) be an irregular minitwistor line in $U$.
Let $C\subset U$ and $\check C\subset\check U$ be as in the proposition. In the following, we use the notations used in the construction of the modification and transformations of the minitwistor lines.
Since $C_0\sim C$ on $U$ as above,
pulling back by the minimal resolution $\mu:\tilde U\lras U$,
we have $l_1 + \ol l_1 +D_{\lmd}+ E + \ol E\sim \tilde C$ on $\tilde U$.
Pulling this back further by the map $\boldsymbol\nu:
\hat U\lras \tilde U$, we have, on $\hat U$, 
\begin{align}\label{}
l_1 + \ol l_1 + D_{\lmd} + E + \ol E + \DDD_0 + \ol \DDD_0
\sim 
\hat C + \DDD + \ol \DDD.
\end{align}
This can be rewritten as
\begin{align}\label{}
\big(l_1 + \ol l_1 + D_{\lmd} + E + \ol E\big) - \hat C
\sim
(\DDD + \ol\DDD)
- (\DDD_0 + \ol \DDD_0).
\end{align}
From this, it is easy to see that 
a linearly equivalence $\check D_{\lmd} \sim \check C$ on the final surface $\check U$ is equivalent to the linearly equivalence 
\begin{align}\label{lineq1}
\DDD + \ol\DDD
\sim \DDD_0 + \ol \DDD_0\quad{\text{on $\hat U$.}}
\end{align}
Note that $\DDD + \DDD$ and $\DDD_0 + \ol \DDD_0$ are obviously deformation equivalent.

We show that \eqref{lineq1} does not hold.
From the definition of the linearly equivalence, \eqref{lineq1} is equivalent to the existence of a meromorphic function $g$ on 
$\hat U$ such that 
\begin{align}\label{lineq2}
(g) = (\DDD_0 + \ol \DDD_0) - (\DDD + \ol\DDD),
\end{align}
where $(g)$ means the principal divisor associated to the function $g$.
This means that the meromorphic function $g$ is of degree one over the rational curves $l_1$ and $\ol l_1$. Namely, $g$ gives a biholomorphic map onto $\CP_1$ from each of these curves.
Also, \eqref{lineq2} means that $g$ does not have a pole on the divisor $E + \ol E +  D_{\lmd} $. 
Hence, $g$ is constant on this divisor.
Since $D_{\lmd}^2 = 0$ on $\hat U$ and $D_{\lmd}\simeq\CP_1$, 
by a theorem of Kodaira \cite{K62} about displacements of a complex submanifold and also by a theorem of Fisher-Grauert \cite{FG65} about the local triviality of the fibration whose fibers are compact non-singular and all biholomorphic, 
there exists an open neighborhood $X$ of $D_{\lmd}$ in $\hat U$ that is biholomorphic to $D_{\lmd}\times \mathbf D$, where $\mathbf D$ is an open disk in $\CC$ .
Since $g$ has to be constant on any fiber of the projection $D_{\lmd}\times\mathbf D\to \mathbf D$ from the holomorphicity on $X$,
and since $E$ and $\ol E$ intersect $D_{\lmd}$ transversely, this implies that the function $g$ is constant on $X$. 
Since $X$ is an open set of $\hat U$ and $\hat U$ is connected, this means that $g$ is constant on the whole of $\hat U$.
This contradicts \eqref{lineq2}.
Therefore, \eqref{lineq2} and hence \eqref{lineq1} also do not hold.
Thus we have proved $\check D_{\lmd} \not\sim \check C$ on  $\check U$.
Since all real small displacements of $\check D_{\lmd}$ in $\check U$ arise from real small displacements of $C$ in $U$ as minitwistor lines, this means that, on $\check U$, only $\mathbb S^1$-invariant minitwistor lines are linearly equivalent to $\check D_{\lmd}$. Therefore, $\dim |\check D_{\lmd}|=1 $.

Finally, to show $\dim|\check C|=0$, it suffices to prove that if $\check C_1$ and $\check C_2$ are the transformations into $\check U$ of distinct regular minitwistor lines $C_1$ and $C_2$ in $U$ respectively, then 
$\check C_1\not\sim\check C_2$. 
To see this, let $\DDD_1$ and $\ol\DDD_1$ be the two pieces of smooth curves  on $\hat U$ that are components of $\boldsymbol\nu\inv(\tilde C_1)$,
and similarly for $\DDD_2$ and $\ol\DDD_2$.
Then as in the same way for obtaining that $\check D_{\lmd}\sim \check C$ 
is equivalent to \eqref{lineq1}, using $C_1\sim C_2$ on $U$, we obtain that $\check C_1\sim\check C_2$ if and only if $\DDD_1 + \ol\DDD_1\sim
\DDD_2 + \ol\DDD_2$. But the latter cannot hold in the same way that \eqref{lineq1} cannot hold.
Hence, $\check C_1\not\sim\check C_2$.
\proofend

\medskip
From the proof, the conclusions in Proposition \ref{prop:le1} about the dimension of the complete linear system hold not only over the whole surface $\check U$ but also over an arbitrary open subset containing the minitwistor line.
On the other hand, the minitwistor space of the de Sitter space is a smooth quadric and the minitwistor lines on it are smooth hyperplane sections. It is easy to prove that for any twistor line $C$ on the quadric, $\dim |C| = 3$ over any open neighborhood of $C$.
Therefore, Proposition \ref{prop:le1} gives the following corollary.

\begin{corollary}\label{c:ndS}
The EW structures on $W\simeq\mathbb S^2\times I$ obtained in Theorem \ref{t:WEW} are not isomorphic to the de Sitter structure.
\end{corollary}

Proposition \ref{prop:le1} also gives the following property of $\check U$ as a complex surface.

\begin{theorem}\label{t:nc}
The complex surface $\check U$ does not admit a compactification as a complex surface.
\end{theorem}

\proof
Let $X$ be a compactification of $\check U$. By the resolution of singularities, we may suppose that $X$ is smooth.
As $\check U$ contains smooth minitwistor lines, $X$ has a smooth rational curve $C$ satisfying $C^2 = 2$.
By adjunction formula and moving $C$ by Kodaira's theorem about a displacement of a compact submanifold, this readily means that $\kappa(X) = -\infty$ for the Kodaira dimension.
Hence, from the classification of compact complex surface, $X$ is either a rational surface, birational to a ruled surface over a curve of positive genus, or a class VII surface.
Any smooth rational curve $C$ satisfying $C^2 = 2$ on a smooth rational surface satisfies $\dim |C| = 3$.
As $X$ has $C$ satisfying $\dim |C|\le 1$ from Proposition \ref{prop:le1}, this means that $X$ cannot be a rational surface. 
Next, any rational curve on a ruled surface over a curve of positive genus has to be a fiber of the ruling, and therefore any rational curve of a surface birational to the ruled surface has non-positive self-intersection number.
Hence, $X$ cannot be birational to a ruled surface over a curve of positive genus.
Finally, any class VII surface has at most a finite number of rational curves but our surface $\check U\subset X$ has a 3-dimensional family of smooth minitwistor lines.
So $X$ cannot be a class VII surface.
This means that the compactification $X$ does not exist.
\proofend

\subsection{Spacelike geodesics}
\label{ss:slg}
Let $S$ be a Segre surface and $W$ the associated EW space obtained in the previous subsection.
In the following, the signature of the EW structure on $W$ is taken to be $(-++)$. 
In this subsection, we prove that the EW structure on $W$ is spacelike Zoll.
The core of the proof will rely on the description of the minitwistor lines given in Section \ref{s:rmtl}.

For any distinct points $x$ and $y$ on $S^{\sigma}_2$,  set
\begin{align}\label{slg}
\mathfrak C_{x,y}:=
\{
H\in W\set x,y\in H
\}.
\end{align}
This set will be a spacelike geodesic in $W$.
We will first show the following.

\begin{proposition}\label{prop:slg1}
The set $\mathfrak C_{x,y}$ in \eqref{slg} is a smooth 1-dimensional submanifold of 
$W$ if it is not an empty set.
\end{proposition}

To show this, we use the deformation theory of compact curves with nodes on a smooth surface. We briefly recall the theory. 
We refer Sernesi's book \cite{SeBook}, Tannenbaum \cite[Lemma (2.2)]{Tan80}, 
\cite[Proposition (1.6)]{Tan82} and Wahl \cite{W74} for the details.

Let $C\subset S_{\reg}$ be a compact nodal curve on a surface $S$ and $p$ a node of $C$.
Let $N'^{\,p}\subset N:=[C]|_C$ be the equisingular normal sheaf of $C$ in $S$ \cite[\S 4.7.1]{SeBook}, with the equisingularity imposed on $p$.
Let $\ms J_C$ be the Jacobian ideal sheaf of the curve $C$. This is a subsheaf of $\ms O_S$ generated by the two partial derivatives of the defining equation of $C$. We then have an expression
\begin{align}\label{N'2}
N'^{\,p} \simeq N\otimes_{\ms O_C} \ms J_C.
\end{align}

The small equisingular displacements of $C$ in $S$ with the equisingularity imposed on $p$ may be investigated by the cohomology groups of the sheaf $N'^{\,p}$.
In particular, if $H^1(N'^{\,p})=0$, then
 writing $d:= \dim H^0(N'^{\,p})$, such displacements are parameterized by a $d$-dimensional complex manifold $Y$ with a reference point $o\in Y$ corresponding to $C$, for which there is a canonical isomorphism $T_oY\simeq H^0(N'^{\,p})$. 
 
When the surface $S$ is compact and rational, which is always the case by \cite[Proposition 2.6]{HN11} for any compact minitwistor space, 
since $H^1(\ms O_S)=0$, 
we readily obtain $\dim H^0(S,[C]) = \dim H^0(N) + 1$.
This means there is a natural isomorphism $T_o |C|\simeq H^0(N)$.
Further, in that case, $Y$ is naturally embedded in the complete linear system $|C|$, so there is an injection $T_oY\subset T_o|C|$.
Under the isomorphisms $T_oY\simeq H^0(N'^{\,p})$ and $T_o |C|\simeq H^0(N)$, 
this injection is equal to the inclusion 
$H^0(N'^{\,p})\subset H^0(N)$ induced from the sheaf inclusion
$N'^{\,p}\subset N$. 

The case where $S$ is the Segre surface and $C$ is a minitwistor line on $S$ is as follows.
 
\begin{proposition}\label{p:h0h1}
Let $S$ be a Segre surface as before and $C\in W$ a regular or an irregular minitwistor line in $S$. Let $p$ be the real node of $C$. Then 
$$
H^1(N'^{\,p}) = 0 \qandq
H^0(N'^{\,p}) \simeq \CC^3.
$$
\end{proposition}

\proof
We follow the argument in \cite[Example 4.7.1 (ii)]{SeBook}.
Let $\nu:\tilde C\lras C$ be the normalization of the nodal curve $C$ at the node $p$,
and write $p'$ and $\ol p'$ for the two points $\nu\inv(p)$ of $\tilde C$.
If $C$ is regular, then $\tilde C\simeq\CP_1$.
If $C$ is irregular, then writing $C = l_i + \ol l_i + D_{\lmd}$
for $i=1,2$ and $\lmd\in I_2$,  $\tilde C$ still consists of three components and we use the same letters for them.
Then for both cases, from \eqref{N'2} and $C^2 = 4$, we have 
\begin{align}
\nu^*N'^{\,p} &\simeq (\nu^*N) \otimes (\nu^*\ms J_C|_C)\notag\\
&\simeq (\nu^*N) \otimes \ms O_{\tilde C}(-p'-\ol p')\label{equi1}.
\end{align}
From this, the projection formula means
\begin{align}
\nu_* \big(\nu^*N'^{\,p}\big) &\simeq N\otimes \nu_*\ms O_{\tilde C}(-p'-\ol p')\notag\\
&\simeq N\otimes \ms J_C|_C\notag\\
&\simeq N'^{\,p}.\label{pf1}
\end{align}
Moreover, since the sheaf $\nu^*N'^{\,p}$ is invertible from \eqref{equi1} and
$\nu$ is the normalization of a node, if $i>0$, then $R^i\nu_*\big(\nu^*N'^{\,p}\big) = 0$ for the higher direct image sheaf.
By the Leray spectral sequence, this means
\begin{align*}
H^i \big(\nu^*N'^{\,p}\big) &\simeq H^i\big(\nu_* (\nu^*N'^{\,p})\big)\\
&\simeq H^i\big(N'^{\,p}\big),\quad \forall i\ge 0.
\end{align*}
From \eqref{equi1}, this implies isomorphisms
\begin{align}\label{fc0}
H^i\big(N'^{\,p}\big) \simeq 
H^i\big((\nu^*N) \otimes \ms O_{\tilde C}(-p'-\ol p')\big),\quad \forall i\ge 0.
\end{align}
The ingredient sheaf $(\nu^*N) \otimes \ms O_{\tilde C}(-p'-\ol p')$ is invertible. 
Recalling that $C\subset\CP_4$ is of degree four,
if $C$ is regular, then the last sheaf is invertible and of degree $4-2=2$.
If $C$ is irregular, then the same sheaf is still invertible, the degree over $l_i,\ol l_i$ and $D_{\lmd}$ being $0,0$ and $2$ respectively.
It readily follows from these that $H^1((\nu^*N) \otimes \ms O_{\tilde C}(-p'-\ol p')) = 0$ and 
$H^0((\nu^*N) \otimes \ms O_{\tilde C}(-p'-\ol p')) \simeq\CC^3$
regardless of whether $C$ is regular or not.
\proofend

\medskip
By the above general deformation theory, 
Proposition \ref{p:h0h1} means that there exists a locally closed, 3-dimensional complex submanifold $Y$ of $|C|\simeq\CP_4$ with a reference point $o\in Y$, which parameterizes all nearby curves in $S$ that are equisingular displacements of $C$, and they satisfy
\begin{align}\label{tan0}
T_o Y\simeq H^0\big(N'^{\,p}\big)\subset T_o |C| = H^0(N).
\end{align}
In a neighborhood of $o$,
our EW space $W$ can be regarded as the real slice of $Y$ with respect to the induced real structure on $Y$.
Gluing these $Y$ that are obtained from all points of $W$, we obtain a complexification $W_{\CC}$ of $W$.
Clearly this is an open subset of the dual variety of $S$ and it is realized in $\CP_4^*$. 
This complexification will be later used.

\medskip\noindent
{\em Proof of Proposition \ref{prop:slg1}.}
From the definition of $\mf C_{x,y}$ as in \eqref{slg}, if $\ol{xy}^*$ denotes the 2-plane in the dual projective space $\CP_4^*$ formed by hyperplanes passing through $x$ and $y$, then $\mf C_{x,y} = W\cap \ol{xy}^{*}$.
Take any member $C\in \mf C_{x,y}$ and let $H\subset\CP_4$ be the hyperplane such that $C = S\cap H$.
The minitwistor line $C$ can be irregular.
In the following, we identify $C$ and $H$.
We show that the complexified tangent space $(T_H W)^{\CC}\simeq\CC^3$ and $T_H(\ol{xy}^{*})\simeq\CC^2$ 
intersect transversely in $T_H(\CP_4^{*})\simeq\CC^4$.

For the identification of the tangent space $T_H(\CP_4^*)$, if $N=[C]|_C$ is the normal bundle of $C$ in $S$ as before, 
there are natural identifications
\begin{align}\label{tan1}
T_H(\CP_4^*) \simeq T_H(N_{H/\CP_4}) \simeq T_C|C| \simeq H^0(N).
\end{align}
where the last one is as noted right before Proposition \ref{p:h0h1}.
Next, from the first isomorphism in \eqref{tan0} and the last isomorphism in \eqref{pf1}, for the complexified tangent space of $T_H W$,
there are natural isomorphisms
\begin{align}\label{tan2}
(T_H W)^{\CC} \simeq H^0(N'^{\,p}) \simeq H^0(N\otimes\ms J_C|_C).
\end{align}
On the other hand, 
\begin{align}\label{tan3}
T_H(\ol{xy}^*) &\simeq \big\{s\in H^0(N_{H/\CP_4})\set 
s(x) = s(y) = 0\big\}\notag\\
&\simeq \big\{s\in H^0(N)\set 
s(x) = s(y) = 0\big\},
\end{align}
where we used the isomorphisms \eqref{tan1} for the latter isomorphism.
Comparing \eqref{tan2} and \eqref{tan3}, since the two points $x$ and $y$ belong to the sphere $S\us_2$ while the node $p$ of $C$ belongs to another sphere $S\us_1$, 
we obtain $T_H(\ol{xy}^*)\not\subset (T_H W)^{\CC}$.
This means that the intersection $T_H(\ol{xy}^*)\cap (T_H W)^{\CC}$ is transveral in $T_H(\CP_4^*)$.
Taking the real slice for $\sigma$, the two submanifolds $W$ and $\ol{xy}^{*\sigma}$  intersect transversely in $\RP_4^*$ at the point $C\in\mf C_{x,y}$.
\proofend

\medskip
In a 3-dimensional Lorenzian space with a 
$(-++)$ signature, the timelike directions at any point are parameterized by an open disk, while the spacelike directions are parameterized by a M\"obius strip.
From this, we readily obtain:
\begin{proposition}\label{prop:slg0}
Any spacelike geodesic on the EW space $W$ is of the form $\mf C_{x,y}$ for some real distinct
points $x,y\in S^{\sigma}_2$.
\end{proposition}

\proof 
Take any point $w\in W$ and let $C$ be the minitwistor line corresponding to $w$.
From the general construction of the geometric structure on the space of minitwistor line \cite[\S 5]{Hi82}, 
if $\tilde C$ is the normalization of $C$ at the real node,
then any geodesic through $w$ is of the form $\mf C_{x,y}$ for some $x,y\in \tilde C$.
Since the geodesic is real, either
$$
\bullet \, x,y\in \tilde C^{\sigma} {\text{ and }} x\neq y,
\quad
\bullet \, x,y\in \tilde C^{\sigma} {\text{ and }} x = y, {\text{ or }} 
\quad
\bullet \, y = \sigma(x) {\text{ and }}  x\neq y.
$$
Among these, the second one gives null geodesics.
The last ones are parameterized by an (open) hemisphere, and therefore from the remark right before the proposition, they give timelike geodesics.
Hence the first one gives a spacelike geodesic.
Further, since $w\in W$, the real locus of $C$ consists of the node $p\in S^{\sigma}_1$ and a circle in $S^{\sigma}_2$, and only the latter gives the real locus of $\tilde C$.
Hence the two points $x$ and $y$ can be identified with points of $S^{\sigma}_2$.
\proofend

\begin{proposition}\label{prop:slg2}
The set $\mathfrak C_{x,y}$ in $W$ is  compact for any distinct points
$x,y\in S^{\sigma}_2$.
\end{proposition}

\proof
As in Section \ref{ss:rmtl}, we consider the subset $S^{*\sigma}_1= \pr_2(  I(S)^{\sigma} \cap \pr_1\inv(S^{\sigma}_1))\subset S^*$. (See \eqref{diagram:double1}--\eqref{diagram:double2} for the notations.) By definition, a hyperplane $H\subset\CP_4$ belongs to 
$S^{*\sigma}_1$ if and only if it belongs to the real pencil 
$(T_p S)^{*\sigma}$ for some point $p\in S^{\sigma}_1$.
Since $I(S)^{\sigma} \cap \pr_1\inv(S^{\sigma}_1)$ is clearly compact, its image $S^{*\sigma}_1$ is also compact.
Hence $S^{*\sigma}_1$ is closed in $\CP_4^*$.
We show the coincidence
\begin{align}\label{Cxy1}
\ol{xy}^*\cap S^{*\sigma}_1 = \mf C_{x,y}.
\end{align}

First, suppose $H\in \mf C_{x,y}$.
Then since $x,y\in H$, we have $H\in\ol{xy}^*$.
Further, since $H\in W$, $H$ is real and the section $S\cap H$ has a singularity at some point $p\in S^{\sigma}_1$.
These mean $H\in S^{*\sigma}_1$.
Hence $H\in \ol{xy}^*\cap S^{*\sigma}_1$. 
So $\mf C_{x,y} \subset \ol{xy}^*\cap S^{*\sigma}_1$.

Next, suppose $H\in \ol{xy}^*\cap S^{*\sigma}_1$.
Since $H\in S^{*\sigma}_1$, there exists a point $p\in S^{\sigma}_1$ such that $H\in (T_p S)^{*\sigma}$.
As $H\in \ol{xy}^*$, we have $x,y\in H$, and 
both belong to $S^{\sigma}_2$ from the choice of $x,y$.
From Remark \ref{rmk:2pts}, this implies $H\in W$. Hence, $H\in \mf C_{x,y}$, so $\ol{xy}^*\cap S^{*\sigma}_1 \subset \mf C_{x,y}.$ Thus \eqref{Cxy1} is shown.

As seen above, $S^{*\sigma}_1$ is a closed subset of $\CP_4^*$.
Also, $\ol{xy}^*$ is a closed subset of $\CP_4^*$.
Hence, from \eqref{Cxy1}, $\mf C_{x,y}$ is a closed subset of 
$\CP_4^*$. Hence $\mf C_{x,y}$ is compact.
\proofend

\medskip
As an immediate consequence of Propositions \ref{prop:slg1}
and \ref{prop:slg2}, we obtain:

\begin{theorem}\label{thm:Zoll}
The Lorenzian EW structure on the 3-manifold $W$, obtained from the Segre surface $S$, 
is spacelike Zoll. Further, it is real analytic.
\end{theorem}

\proof
By Proposition \ref{prop:slg0}, any spacelike geodesic on $W$ is 
of the form $\mf C_{x,y}$ for some distinct points
$x,y\in S^{\sigma}_2$. 
By Propositions \ref{prop:slg1} and \ref{prop:slg2},  $\mf C_{x,y}$ is always a smooth closed curve in $W$ if it is not empty.
This means that $W$ is spacelike Zoll.

The real analyticity of the EW structure on $W$ is obvious on the locus of regular minitwistor lines for the same reason for real analyticity of the EW structure on the space of smooth minitwistor lines. Real analyticity on the locus of irregular minitwistor lines follows for the same reason because it is obtained as those on the modification $\check U$ whose minitwistor lines are all smooth.
\proofend

\begin{remark}{\em 
As shown in the proof of Proposition \ref{prop:slg1}, the intersection of $W$ and the 2-plane $(\ol{xy})^{*\sigma}$ is transversal.
This implies that none of the spacelike geodesics $\mf C_{x,y}$ are `non-reduced'; in other words, they are not multiply covered by nearby spacelike geodesics.
}
\end{remark}

\section{Moduli, automorphisms, and deformation to de Sitter structure}
\label{s:mad}
\subsection{The family of disks in the Segre surface}\label{ss:disk}
As before, let $W$ be the EW space obtained from the Segre surface $S$, $W_0\subset W$ the dense open subset that parameterizes regular minitwistor lines, and $\ol W= W\cup\ms H_2\cup\ms H_4$ the natural compactification of $W$ (see \eqref{W}).
Let 
\begin{equation}\label{tf1}
   \vcenter{
   \xymatrix@=18pt{
\ms C_0
\ar @{^{(}->}[r]
\ar[d] 
& \ms C
\ar[d]
\ar @{^{(}->}[r] 
& \ol{\ms C}
\ar[d] \\
W_0
\ar @{^{(}->}[r] 
& W
\ar @{^{(}->}[r] 
& \ol W
}
}
\end{equation}
be the tautological families over $W,W_0$ and $\ol W$, so that 
$\ms C, \ms C_0$ and ${\ol{\ms C}}$ are subsets of $S \times W, S\times W_0$ and $S\times \ol W$ respectively and
the fiber over a point $C\in W$ is exactly the curve $C\subset S$, and so on.
For each curve $C\in W$, we write $C^{\RR}$ for the real circle in $C$. 
Using the sphere $S\us_2$, $C^{\RR}=C\cap S\us_2$. 
These circles constitute the loci $\ms C^{\RR}\subset\ms C$ and $\ms C^{\RR}_0\subset\ms C_0$, whose projection to $W$ and $W_0$ gives the structure of $S^1$-bundles.
If $C\in \ol W\minus W=\ms H_2\cup\ms H_4$, then $C$ does not have a real circle anymore, but it has a real point that can be regarded as a shrinking limit of the real circles.
So this real point belongs to $S\us_2$.
There is another real point on $C$ that can be regarded as a limit of the real node of the minitwistor lines, and it belongs to $S\us_1$.
Over the locus $\ol W\minus W$, adding these points to $\ms C^{\RR}$,
we obtain the subset ${\ol{\ms C}}^{\RR}$ of $\ol{\ms C}$ and this is a natural compactification of $\ms C^{\RR}$.

Taking the normalization for all fibers of \eqref{tf1} at the real node belonging to $S\us_1$, we obtain 
the families $\ms C'\lras W,\, \ms C'_0\lras W_0$ and 
$\ol{\ms C}'\lras \ol W$ whose fibers are the normalizations.
Since $S\us_1\cap S\us_2=\emptyset$, 
the loci $\ms C^{\RR},\,\ms C^{\RR}_0$ and ${\ol{\ms C}}^{\RR}$ can be thought as subsets of ${\ms C}',\, {\ms C}'_0$ and $\ol{\ms C}'$ respectively.
The structure of all fibers of the projections
from the complements ${\ms C}'\minus \ms C^{\RR},\, {\ms C}'_0\minus \ms C^{\RR}_0$ and $\ol{\ms C}'\minus{\ol{\ms C}}^{\RR}$ can be readily identified.
In particular, all fibers consist of two connected components. For example, those of the projection $\ms C'_0\minus \ms C^{\RR}_0\lras W_0$ consist of two disjoint open disks.
As $W\simeq\mathbb S^2\times I$ is simply connected, 
these three complements consist of two connected components. 
Attaching back $\ms C^{\RR},\,\ms C^{\RR}_0$ and ${\ol{\ms C}}^{\RR}$
to each of these connected components independently, we obtain relative compactifications of these components.
Let $\ms C^+$ and $\ms C^-$ be such compactifications of the two components of $\ms C'\minus\ms C^{\RR}$, and similarly for 
${\ms C}'_0\minus \ms C^{\RR}_0$ and $\ol{\ms C}'\minus{\ol{\ms C}}^{\RR}$.
The fibers of the projections $\ms C^{\pm}_0\lras W$ are closed disks, while the fibers of $\ol{\ms C}^{\pm}\lras \ol W$ over a point of $\ms H_j$ ($j=2,4$)  
can be written $\pi_j\inv(\CCC)$ or $\pi_j\inv(\ol\CCC)$ where $\CCC$ and $\ol\CCC$ are lines on $Q_j$ intersecting at a point as in \eqref{AolA} because the decomposition \eqref{AolA} can be regarded as a limit of the decomposition of a regular minitwitor line by the real circle.

The connected components $\ms C^+_0$ and $\ms C^-_0$ can be characterized by the following intersection property with some lines on $S$.

\begin{lemma}\label{l:+-} 
Exactly one of the following two situations occurs.
(i) Any $C^+\in \ol{\ms C}^+_0$ intersects the lines $l_1$ and $\ol l_2$ and does not intersect $\ol l_1$ and $l_2$.
(ii) Any $C^+\in \ol{\ms C}^+_0$ intersects the lines $\ol l_1$ and $l_2$ and does not intersect $l_1$ and $\ol l_2$. 
\end{lemma}

Note that this property cannot hold for ${\ms C}^+$ nor $\ol{\ms C}^+$ because any irregular minitwistor line contains either $l_1+\ol l_1$ or $l_2+\ol l_2$.

\medskip\noindent
{\em Proof of Lemma \ref{l:+-}.}
We put $\ol W_0:= \ol W\minus(\ol\varpi\inv(p_1)\cup\ol\varpi\inv(p_2))$.
Let $i\in\{1,2\}$.
Since any $C\in \ol W_0$ is a hyperplane section of the Segre surface $S$ and contains none of $l_i$ nor $\ol l_i$ $(i=1,2)$,
both $C\cap l_i$ and $C\cap \ol l_i$ consist of one point, and they are not equal since $p_i\not\in C$.
Hence, exactly one of $C^+\cap l_i\neq\emptyset$ and $C^-\cap l_i\neq\emptyset$ holds.
This means that, if we put
$$
\ol W_{0,i}^+=\{C\in \ol W_0\set C^+\cap l_i \neq \emptyset\}
\qandq
\ol W_{0,i}^-=\{C\in \ol W_0\set C^-\cap l_i \neq \emptyset\},
$$
then 
$\ol W_0 = \ol W_{0,i}^+\sqcup \ol W_{0,i}^-$ ($i=1,2$).
Further, if $C\in \ol W_{0,i}^+$, then since the point $C^+\cap l_i$ belongs to the interior of the closed disk $C^+$ and is not $p_i=l_i\cap \ol l_i$, $\ol W_{0,i}^+$ is an open subset of $\ol W_0$.
Similarly, $\ol W_{0,i}^-$ is also open in $\ol W_0$.
Hence, from the connectedness of $\ol W_0$, either $\ol W_0 = \ol W_{0,i}^+$ or $\ol W_0 = \ol W_{0,i}^-$ for $i=1,2$.
So there are four possibilities.

To complete a proof, it suffices to show that 
either $\ol W_{0,1}^+\cap \ol W_{0,2}^-\neq\emptyset$
or $\ol W_{0,1}^-\cap \ol W_{0,2}^+\neq\emptyset$ because these respectively mean (i) and (ii) in the lemma.
We show this by making use of the boundaries $\ms H_j$.
Let $j\in\{2,4\}$ and take any $C=\pi_j\inv(\CCC)+\pi_j\inv(\ol\CCC)\in \ms H_j$ as in \eqref{AolA} that does not belong to $\ol\varpi\inv\{p_1,p_2\}$.
Then from what we have remarked right before the present lemma, either $\pi_j\inv(\CCC)\in \ol{\ms C}^+_0$ or $\pi_j\inv(\CCC)\in \ol{\ms C}^-_0$.
Since the bidegrees of the lines $\pi_j(l_1)$ and $\pi_j(\ol l_2)$ in the quadric $Q_j\simeq\qdr$ are equal from the explicit description in Section \ref{s:S}, either ``$\pi_j\inv(\CCC)\cap l_1\neq\emptyset$
and $\pi_j\inv(\CCC)\cap \ol l_2\neq\emptyset$''
or ``$\pi_j\inv(\CCC)\cap l_1=\emptyset$
and $\pi_j\inv(\CCC)\cap \ol l_2=\emptyset$''.
If $\pi_j\inv(\CCC)\in \ol{\ms C}^+_0$ and $\pi_j\inv(\CCC)\cap l_1\neq\emptyset$, then $C\in \ol W_{0,1}^+\cap \ol W_{0,2}^-$, and 
if $\pi_j\inv(\CCC)\in \ol{\ms C}^+_0$ and 
$\pi_j\inv(\CCC)\cap l_1=\emptyset$, then 
$C\in \ol W_{0,1}^-\cap \ol W_{0,2}^+$.
The remaining case  $\pi_j\inv(\CCC)\in \ol{\ms C}^-_0$ can be shown in the same way.
\proofend

\begin{definition}\label{d:pm}{\em
In the following, we distinguish between the two families $\ms C^+_0$ and $\ms C^-_0$ of disks, where {\em members of $\ms C^+_0$ intersect the line $\ol l_1$ (and hence the line $l_2$ as well)}.
This distinction is well-defined by Lemma \ref{l:+-}, and naturally extends to those of $\ms C^+$ and $\ms C^-$ as well as those of $\ol{\ms C}^+$ and $\ol{\ms C}^-$. } 
\end{definition}

From the previous lemma,
we have the inclusion $\cup_{C\in W_0} C^+\subset S\minus (l_1\cup \ol l_2)$.
The next lemma, which is also important for our purpose, means reverse inclusion.

\begin{lemma}\label{l:disk1}
For any point $z\in S\minus (l_1\cup \ol l_2)$, there exists a point $C\in W_0$ such that $z\in C^+$. 
\end{lemma}

We show this lemma in two steps.
We first prove the case $z\in S\us$, which is much easier. Since $z\not\in l_1\cup \ol l_2$, this implies $z\not\in\{p_1,p_2\}$.

\medskip
\noindent{\em Proof of Lemma \ref{l:disk1} in case $z\in S\us$.}
Recall that $S\us=S\us_1\sqcup S\us_2$.
If $z\in S\us_1$, then noting that $z\neq p_1, p_2$ as above, the interval $\varpi\inv(z)=J_4(z)$ is included in $W_0$.
Let $C$ be any point of the interval $J_4(z)$. The real node of the minitwistor line $C$ is exactly $z$. 
Hence, $z\in C^+\cap C^-$, so $z\in C^+$. 
Next, to show the case $z\in S\us_2$, recall from Section \ref{s:rmtl} that for any $p\in S\us_1$, the family of circles $\{ S\us_2\cap C\set C\in J_4(p)\}$
almost foliates the sphere $S\us_2$ in the sense  
\begin{align}\label{union22}
\bigsqcup_{C\in J_4(p) } S^{\sigma}_2\cap C = S^{\sigma}_2\minus\{\tau_2(p),\tau_4(p)\}.
\end{align}
(See \eqref{union2} and Figure \ref{fig:Ssigma2}.)
Since the disks $C^+$ and $C^-$ contain the boundary circle, it follows from \eqref{union22}
that any point $z$ of $S^{\sigma}_2\minus\{\tau_2(p),\tau_4(p)\}$ is passed by $C^+$ for some $C\in J_4(p)$. 
The missing points $\tau_2(p)$ and $\tau_4(p)$ are known to be passed through 
by $(C')^+$ for some $C'\in  J_4(p')$ 
just by taking another point $p'\in S\us_1$, using \eqref{union22} for the point $p'$ and $\{\tau_2(p),\tau_4(p)\}\cap \{\tau_2(p'),\tau_4(p')\}=\emptyset$.
\proofend

\medskip
To prove the case $z\not\in S\us$, we recall from Section \ref{s:S} that for $j=2,4$ and $1\le i\le 4$, the images $\pi_j(l_i)$ and $\pi_j(\ol l_i)$ are $\CC^*$-invariant lines on $Q_j$ and they form one distinguished ``square'' in $Q_j$.
The ``vertices'' of the square are exactly the $\CC^*$-fixed points on $Q_j$.
In the following argument, we use the words `square,' `edge' and `vertex' to mean these. 
If we further assume $z\not\in (l_1\cup \ol l_2)$ as in the present situation, then 
$z\neq e_0,e_1$ and this means that the point $\pi_j(z)$ ($j=2,4$) is not the vertex of the square.
This implies that there exists at least one line $\CCC$ on $Q_j$ through the non-real point $\pi_j(z)$, which is not an edge of the square.
Put $C_j:=\pi_j\inv(\CCC)$ and $\ol C_j:=\pi_j\inv(\ol\CCC)$, both of which are irreducible conics as $\CCC$ and $\ol \CCC$ are not edges.
Obviously, $z\in C_j$ and $z\not\in \ol C_j$.
The curve $C: = C_j + \ol C_j$ belongs to the boudnary $\ms H_j$
 (see \eqref{AolA}).
Exchanging $\CCC$ and $\ol \CCC$ if necessary, we may suppose $C_j\in \ol{\ms C}^+_0$. Then $z\in C^+$.
The intersection $C_j\cap \ol C_j$ consists of two points, one of which belongs to $S\us_1$ and the other to $S\us_2$. 
If $p$ means the former point, then $\tau_j(p)$ is the latter point.
We shall show Lemma \ref{l:disk1} in case $z\not\in S\us$ by proving that the curve $C=C_j+\ol C_j\in \ms H_j$ ($j=2,4$) admits a displacement in $S$ preserving the real structure that remains to pass through $z$ and that is equisingular only at the node $p$.

\medskip
\noindent{\em Proof of Lemma \ref{l:disk1} in case $z\not\in S\us$.}
We use the notations in Proposition \ref{prop:slg1}.
In particular, $\nu:\tilde C\lras C$ still means the normalization of $C$ at the node $p$, so that $\tilde C$ still has another node $\tau_j(p)$.
The above equi-singular displacements of $C$ may be investigated through the cohomology groups $H^i(N'^{\,p})$, and if $p'\in C_j\subset \tilde C$ and $\ol p'\in \ol C_j\subset \tilde C$ mean the two points $\nu\inv(p)$, then in the same way to \eqref{fc0}, we obtain isomorphisms
$$
H^i\big(N'^{\,p}\big )\simeq H^i\big((\nu^*N)\otimes \ms O_{\tilde C}(-p'-\ol p')\big)\quad\forall i\ge 0.
$$
The sheaf $(\nu^*N)\otimes \ms O_{\tilde C}(-p'-\ol p')$ is invertible and the degrees over $C_j$ and $\ol C_j$ are both one.
From this, using the normalization of $\tilde C$ at the remaining node $\tau_j(p)$, we readily obtain 
$$
H^i\big((\nu^*N)\otimes \ms O_{\tilde C}(-p'-\ol p')\big)
\simeq
\begin{cases}
\CC^3 & i=0,\\
0 & i>0.
\end{cases}
$$
Hence, the versal family of the above equi-singular displacements of $C$ in $S$ is parameterized by a 3-dimensional complex manifold whose tangent space at the point corresponding to $C$ is identified with $H^0((\nu^*N)\otimes \ms O_{\tilde C}(-p'-\ol p'))$.

Next, we consider equi-singular displacements of the same curve $C$ in $S$ where this time equi-singularity is imposed at both nodes $p$ and $q:=\tau_j(p)$ so that the displacements remain reducible.
In a similar way to the above,
if $N'^{\,p,q}$ denotes the equisingular normal sheaf of $C$ in $S$ at $p$ and $q$,
then these equi-singular displacements of $C$ in $S$ are investigated through $H^i(N'^{\,p,q})$.
If $\nu: \tilde C \lras C$ is the full normalization of $C$ this time and $q'$ and $\ol q'$ are the points of $C_j$ and $\ol C_j$ (regarded as components of $\tilde C$) respectively that are mapped to $q$ by $\nu$, then using the points $p'\in C_j\subset \tilde C$ and $\ol p'\in \ol C_j\subset \tilde C$,
\begin{align}\label{N'ds}
\nu^* N'^{\,p,q} \simeq
(\nu^*N)|_{C_j}
\otimes\ms O_{C_j}(-p'-q')
\,\oplus\,
(\nu^*N)|_{\ol C_j}
\otimes\ms O_{\ol C_j}(-\ol p'-\ol q').
\end{align}
As before, the two direct summands in the right-hand side of \eqref{N'ds}
are invertible, and further, as $C_j$ and $\ol C_j$ are conics,
they are of degree zero (i.e.\,trivial).
From this, using the Leray spectral sequence, we obtain $H^1(N'^{\,p,q}) =0$ 
and $H^0(N'^{\,p,q}) \simeq\CC^2$.
Hence, the versal family of the present equi-singular displacements is parameterized by a non-singular complex surface whose tangent space at the point $C$ is 
identified with $H^0(N'^{\,p,q})$.
Since the boundary $\ms H_j$ indeed parameterizes such equi-singular displacements that preserve the real structure,
a neighborhood of the point $C$ in the boundary $\ms H_j$ is a real slice of the versal family.
So we write $\ms H_j^{\CC}$ for the parameter space of this versal family around the point $C$.

Next, we consider equi-singular displacements of $C$ in $S$ with respect to the node $p$ which remain to pass through the two points $z\in C_j$ and $\ol z\in \ol C_j$.
The relevant sheaf is $N'^{\,p}\otimes\ms O_C(-z-\ol z)$.
Recall that $\{z,\ol z\}\cap \{p,q\}=\emptyset$.
Pulling the sheaf back by the above full normalization $\nu$,
we obtain 
$$
\nu^*\big(N'^{\,p}\otimes\ms O_C(-z-\ol z)\big)
\simeq
(\nu^*N)|_{C_j}\otimes \ms O_{C_j}(-p'-z)
\,\oplus\,
(\nu^*N)|_{\ol C_j}\otimes \ms O_{\ol C_j}(-\ol p'-\ol z).
$$
Again, the two direct summands on the right-hand side are trivial.
But the cohomology groups of 
$N'^{\,p}\otimes\ms O_C(-z-\ol z)$ 
is not isomorphic to those of the pullbacks by $\nu$, and instead, we have to identify the stalks over $q'$ and $\ol q'$ to get a sheaf over the partial normalization of $C$ at $p$ and to consider the cohomology groups of this sheaf, which turns out to be isomorphic to the cohomology groups of $N'^{\,p}\otimes\ms O_C(-z-\ol z)$.
From this, using $\{z,\ol z\}\cap \{p,q\}=\emptyset$, we obtain $H^1(N'^{\,p}\otimes\ms O_C(-z-\ol z))=0$ and further, 
$H^0(N'^{\,p}\otimes\ms O_C(-z-\ol z))\simeq\CC$.
Moreover, regarding both $H^0(N'^{\,p}\otimes\ms O_C(-z-\ol z))$ and 
$H^0(N'^{\,p,q})$ as subspaces of $H^0(N'^{\,p})\simeq\CC^3$,
by the natural inclusions, we have
$$
H^0\big(N'^{\,p}\otimes\ms O_C(-z-\ol z)\big)\cap
H^0\big(N'^{\,p,q}\big) = 0.
$$
This means that the parameter spaces of the two versal families intersect transversely at the point $C$ in the 3-dimensional parameter space of the versal family of equi-singular displacements of $C$ in $S$ with respect to the node $p$.
Let $\mf C^{\CC}$ be the parameter space of the versal family of equi-singular displacements of $C$ in $S$ with respect to the node $p$
fixing the two points $z$ and $\ol z$, so that $\mf C^{\CC}$ is a smooth holomorphic curve in the last 3-dimensional parameter space through the reference point corresponding to $C$.
As above, this curve and the complexification $\ms H_j^{\CC}$ intersect transversely at the point $C$.

Since the real structure on $S$ induces that on the sheaf $N'^{\,p}\otimes\ms O_C(-z-\ol z)$, the parameter space $\mf C^{\CC}$ is also equipped with a real structure.
Let $\mf C$ be the real locus of this real structure.
Since $\mf C^{\CC}$ is smooth, $\mf C$ is necessarily a (real) smooth curve.
By the above transversality, just by taking the real locus, $\mf C$ intersects the boundary $\ms H_j$ transversely at the point $C$.
From the results in Section \ref{ss:var1}, 
for any real curve $C'\subset S$ belonging to a neighborhood in the real locus $S^{*\sigma}_1$ of the point $C\in \ms H_j$,
$C'$ has a real circle (in $S\us_2$) if and only if $C'$ belongs to the open 3-manifold $W$.
Hence, redefining $\mf C$ as the part of the above $\mf C$ that is contained in $W$, 
$\mf C$ is a smooth curve in $W$ whose end is the point $C\in \ms H_j$.
Any $C'\in \mf C$ satisfies $z\in C'$.
Further, since $z\in C_j\in\ol{\ms C}^+_0$, we have $(C')^+\in \ms C^+_0$ for $C'\in \mf C$. 
This completes the proof of Lemma \ref{l:disk1} in the case $z\not\in S\us$.
\proofend

\begin{remark}
{\em
The curve $\mf C$ in the proof is a timelike geodesic in $W$ emanating from future or past infinity. (See the proof of Proposition \ref{prop:slg0}.)
}\end{remark}

\medskip 
In the rest of this paper, we denote 
\begin{align}\label{Splus}
S^+:= S\minus(l_1\cup \ol l_2).
\end{align}
This is an open smooth complex surface but it does not have a real structure.
From Lemmas \ref{l:+-} and \ref{l:disk1}, we obtain the following corollary, which is important for us.

\begin{corollary}\label{c:disk} 
Any diffeomorphism of $W_0$ preserving the conformal structure induces a diffeomorphism of $S^+$.
\end{corollary}

In the next subsection, we will show that this diffeomorphism of $S^+$ is holomorphic if the diffeomorphism of $W_0$ preserves the EW structure.

\medskip
\noindent{\em Proof of Corollary \ref{c:disk}.}
From Lemmas \ref{l:+-} and \ref{l:disk1}, we have
$$\bigcup_{C\in W_0} C^+ = S^+.$$
This means that the complex surface $S^+$ is the minitwistor space of the EW space $W_0$ whose minitwistor lines are closed disks. 
Of course, $S^-:=S\minus(\ol l_1\cup l_2)$ is the minitwistor space whose minitwistor lines are other disks $C^-$.
Hence, if $\phi:W_0\lras W_0$ is any diffeomorphism preserving the conformal structure, then the differential of $\phi$ induces either a diffeomorphism $S^+\lras S^+$ or $S^+\lras S^-$. If the latter is the case, then by taking the composition with the real structure $\sigma$ on $S$, since $\sigma$ flips $C^+$ and $C^-$, we obtain a diffeomorphism $S^+\lras S^+$.
\proofend

\medskip
Since the family $\ms C^+_0$ is realized in the product $S^+\times W_0$, 
there is a double fibration
\begin{align}\label{d:double3}
 \xymatrix{ 
&\ms C^+_0 \ar[dl]_{\pr_1} \ar[dr]^{\pr_2} &\\
S^+ && W_0
 }
\end{align}
All fibers of $\pr_2$ are (smooth) closed disks, which are the $+$-part of the regular minitwistor lines.
The Einstein-Weyl structure $([g],\nabla)$ is defined so that 
\begin{itemize}
\item[(i)]  $\pr_2(\pr_1^{-1}(q))$ is a timelike geodesic if $q\in S^+\setminus S^\sigma_2$,
\item[(ii)] $\pr_2(\pr_1^{-1}(q))$ is a null surface if $q\in S^\sigma_2$.
\end{itemize}

Here we recall some properties concerning the correspondence \eqref{d:double3}. 
In general, on the 3-dimensional torsion-free indefinite Weyl space $(X,[g], \nabla)$,   
a tangent plane $\Pi^\RR \subset T_x X$ at $x\in X$ is called {\it null} if 
$[g]$ degenerates on $\Pi^\RR$. A real surface $\mathcal S \subset X$ is called {\it null} 
if its tangent space $T_x \mathcal S \subset T_xX$ is null for any $x\in \mathcal S$. 

On the correspondence \eqref{d:double3}, 
the Lorenzian conformal structure $[g]$ on the space $W_0$ is defined so that the condition (ii) above holds. 
The torsion-free connection $\nabla$ is defined so that each null surface is totally geodesic. 
Further, condition (ii) guarantees the following integrability condition:
for any null plane $\Pi^\RR$ at any point $C\in W_0$, 
there exists a null surface $\mathcal S$ containing $C$ such that $T_C \mathcal S= T_C\Pi^\RR $. 
This integrability condition is equivalent to the Einstein-Weyl condition. 
By construction, any maximal null surface on $W_0$ is written as $\pr_2(\pr_1^{-1}(q))$
for some $q\in S^\sigma_2$.

\subsection{The inverse correspondence}\label{ss:inv}
We show that any automorphism $\phi: W_0 \to W_0$ preserving the EW structure 
induces a holomorphic isomorphism $\phi^+ : S^+ \to S^+$.
This fact means that the complex structure on $S^+$ is characterized by the EW structure 
$([g], \nabla)$ on $W_0$, that is, the inverse correspondence works.

Let us fix a real Lorenzian metric $g\in [g]$ on $TW_0$ with signature $(-++)$. 
We extend $g$ to a complex bi-linear form on $T^\CC W_0$. 
We define 
$$ \mathcal N= \{ [v] \in \PP(T^\CC W_0) \mid g(v,v)=0\}$$ 
and let $p_2 : \mathcal N \to W_0$ be the projection. 
For each $x = [v] \in \mathcal N$, we write $\bar x = [ \bar v ] \in \mathcal N$ where  
$\bar v$ is the complex conjugation of $v$ in $T^\CC W_0$. 
The fixed point set of this conjugation map is 
$$ \mathcal N^\RR = \{ [v] \in \PP(TW_0) \mid g(v,v)=0\}. $$
The natural projection $p_2 : (\mathcal N, \mathcal N^\RR)\to W_0$ is a $(\CP_1, \RP_1)$-bundle.
Then  $\mathcal N$ is divided into two closed disk bundles $\mathcal N^+$ and $\mathcal N^-$ 
so that each fiber of $\mathcal N^\RR \subset \mathcal N^\pm$ is the boundary circle  and 
$\mathcal N^+ \cap \mathcal N^- = \mathcal N^\RR$. 
Here $\mathcal N^+$ and $\mathcal N^-$ are exchanged by the conjugation. 
We show that $\mathcal N^+$ is naturally identified with $\mathscr C^+_0$. 
For this, we notice the following elementary properties. 

\begin{lemma}  \label{lem_LA}
Let $g$ be a real Lorenzian metric on $\RR^3$ with signature $(-++)$. 
\begin{itemize}
\item If $v\in \RR^3$ is non-zero and null, then 
$\Pi_v^\RR = \{ w \in \RR^3 \mid g(w,v)=0 \}$ is the unique null plane containing $v$. 
Conversely, any null plane contains a unique null direction. 
\item If $v\in \RR^3$ is spacelike, i.e. $g(v,v)>0$, there exists exactly two null planes 
containing $v$. Conversely, any two distinct null planes intersect along a spacelike line.  
\item If $v\in \RR^3$ is timelike, i.e. $g(v,v)<0$, there are no null planes containing $v$. 
If we consider the complex bi-linear extension of $g$ on $\CC^3 = \RR^3 \otimes \CC$, 
there are just two complex null planes containing $v$, and these planes are 
complex conjugate to each other. Conversely, let $\Pi \subset \CC^3$ be any complex null plane 
such that $\Pi \neq \overline \Pi$, then $l^\RR = \Pi \cap \RR^3$ is a real timelike line. 
We also obtain $l^\RR \otimes \CC = \Pi \cap \overline \Pi$. 
\end{itemize}
\end{lemma}

We define a natural map $p_1: \mathcal N^+ \to S^+$ in the following way. 
For any $x =[v] \in \mathcal N^+$, 
$x$ is a complex null line on $T^\CC_C W_0$ where $C= p_2(x) \in W_0$. 
Let 
$$ \Pi_{x} = \{ w \in T^\CC_C W_0 \mid g(w,v)=0 \}$$
be the unique complex null plane containing $x$. 

If $x =[v] \in \mathcal N^\RR$, we can assume $v$ is real, i.e. $v\in T_C W_0$.
Then, we can also consider the real null plane
$$ \Pi_x^\RR = \{ w \in T_C W_0 \mid g(w,v)=0 \}. $$ 
Notice that $\Pi_x  = \Pi_{\overline{x}} = \overline { \Pi_{x}}$
and $\Pi_{x} = \Pi_{x}^\RR \otimes \CC$. 
By the integrability of Einstein-Weyl structure, 
there exists a maximal null surface $\mathcal S \subset W_0$, 
such that $T_C\mathcal S = \Pi_x^\RR$. 
By the correspondence \eqref{d:double3}, we can write 
$\mathcal S = \pr_2(\pr_1^{-1}(q))$ for a uneque $q\in S^\sigma_2$. 
Then we define $p_1(x)=q$. 
Notice that for any null surface $\mathcal S \subset W_0$, 
we can define its natural lift $\widetilde {\mathcal S} \subset \mathcal N^\RR$ 
so that each point $x\in \widetilde {\mathcal S}$ correspondes with the 
null plane $T_{p_2(x)}\mathcal S \subset T_{p_2(x)} W_0$. 
Hence $\mathcal N^\RR$ is foliated by such real surfaces, 
and the map $p_1$ is the natural projection of this foliation to the leaf space $S^\sigma_2$. 

If $x\in \mathcal N^+\setminus \mathcal N^\RR$, then $\Pi_x \neq \overline{\Pi_{x}}$. 
This is the third case of Lemma \ref{lem_LA}, 
so $l_x = \Pi_x \cap \overline{\Pi_{x}}$ is a complex line, 
$l_x^\RR =\Pi_x \cap T_CW_0$ is a real timelike line, and we obtain
 $l_x = l_x^\RR \otimes \CC$.
Hence each $x\in \mathcal N^+ \setminus \mathcal N^\RR$ one-to-one corresponds with 
a timelike line on $W_0$. 
For any timelike geodesic $\mathcal C \subset W_0$, we obtain its natural lift 
$\widetilde{\mathcal C} \subset \mathcal N^+ \setminus \mathcal N^\RR$. 
So $\mathcal N^+ \setminus \mathcal N^\RR$ is foliated by such lifts, 
and its leaf space is the set of timelike geodesics on $W_0$, that is $S^+ \setminus S^\sigma_2$. 
So we obtain a natural projection $p_1: \mathcal N^+ \setminus \mathcal N^\RR \to S^+ \setminus S^\sigma_2$. 
Identifying $\mathcal N^+$ with the family $\mathscr C^+_0$ naturally,
from \eqref{d:double3}, we obtain the double fibration

\begin{align}\label{d:double4}
 \xymatrix{ 
&(\mathcal N^+, \mathcal N^\RR) \ar[dl]_{p_1} \ar[dr]^{p_2} &\\
S^+ && W_0
 }
\end{align}

From Corollary \ref{c:disk}, any diffeomorphism $\phi: W_0 \to W_0$ preserving the conformal structure 
induces a diffeomorphism $\phi^+:S^+\lras S^+$.
To show that $\phi^+$ is holomorphic when $\phi$ preserves the Einstein-Weyl structure, 
we study the correspondence more precisely in terms of distributions. 
We define a complex rank two distribution $E$ on $\mathcal N$ 
so that $E_x$ is the horizontal lift of the null plane $\Pi_x$ at $x\in \mathcal N$. 
We also define a complex rank 3 distribution $F = E\oplus V^{0,1}$ on $\mathcal N$
where $V^{0,1}$ is the bundle of $(0,1)$-tangent vector of the fiber of the $\CP_1$-bundle $p_2: \mathcal N \to W_0$.  
By the general theory, the integrability of $F$ is equivalent to the EW condition (see \cite{Hi82} or \cite{Nkt09}). 
We show the following
\begin{proposition} \label{prop:01t}
For each point $x\in \mathcal N^+ \setminus \mathcal N^\RR$, 
the image $ (p_1)_* (F_x)$ is equal to the $(0,1)$-tangent space $T^{0,1}_{p_1(x)}S^+$. 
\end{proposition}

This readily means:

\begin{corollary}\label{c:01t}
If $\phi$ is a diffeomorphism of $W_0$ preserving the EW structure, 
then the induced diffeomorphism $\phi^+:S^+\lras S^+$ is holomorphic.
\end{corollary}

\proof
$\phi$ induces an automorphism of $(\mathcal N, \mathcal N^\RR)$ preserving the distribution $E$ and $F$. 
Hence the induced diffeomorphism $\phi^+: S^+ \to S^+$ preserves $(0,1)$-tangent bundle $T^{0,1}S^+$ 
by Proposition \ref{prop:01t}.
Thus $\phi^+$ is holomorphic. 
\proofend

\medskip
To prove  Proposition \ref{prop:01t}, 
let us define the real rank 1 distribution $L^\RR$ on $\mathcal N^+ \setminus \mathcal N^\RR$,
so that $L^\RR_x$ is the horizontal lift of the timelike line $l^\RR_x$ at $x$. 
We define $L= L^\RR \otimes \CC$, then we obtain 
$$ L\subset E \subset F, \qquad L = E \cap \overline E = F \cap \overline F. $$
By the condition (i) of the double fibration \eqref{d:double3}, we obtain 
$ L = \ker(({p_1})_* \otimes \CC). $

\begin{proof}[Proof of Proposition \ref{prop:01t}]
It is obvious that 
$(p_1)_*(V^{0,1}_x) \subset T^{0,1}_{p_1(x)} S^+$ for any $x\in \mathcal N^+\setminus \mathcal N^\RR$ 
since the restriction of $p_1$ to each fiber of $p_2$ is a rational curve. 
So it is enough to show 
\begin{equation} \label{eq:01t}
(p_1)_* (E_x) \subset T^{0,1}_{p_1(x)} S^+. 
\end{equation}
Actually, if \eqref{eq:01t} holds, 
we obtain $(p_1)_*(F_x) = (p_1)_*(E_x\oplus V^{0,1}_x) \subset T^{0,1}_{p_1(x)} S^+$. 
The statement follows since   
$\dim_\CC (p_1)_*(F_x) = \dim_\CC (F_x/L_x) = 2 = \dim_\CC T^{0,1}_{p_1(x)}S^+$.

To show \eqref{eq:01t}, we study the complex minitwistor correspondence. 
By a similar construction to the real case, we obtain the following diagram. 
\begin{align}\label{d:double6}
 \xymatrix{ 
&\mathcal N^+ \ar[dl]_{p_1} \ar[dr] \ar[r]^\iota & \mathcal N^\CC \ar[dr]^{\tilde p_2} \ar[dll]^{\tilde p_1} & \\
S^+ && W_0 \ar[r]^{i} & W_0^\CC
 }
\end{align}
Here $W^{\CC}_0$ is a complex 3-fold with complex Einstein-Weyl structure $([g^\CC],\nabla)$, 
and $\tilde p_2 : \mathcal N^\CC \to W^\CC_0$ is the $\CP_1$-bundle defined by
\begin{align}
 \mathcal N^\CC = \{ [v] \in \PP(T^{1,0}W_0^\CC) \mid g^\CC(v,v)=0 \}.
\end{align}
The EW manifold $W^\CC_0$ is constructed as the set of 
complex regular minitwistor lines, that is, nodal rational curves 
on $S^+$ such that they have a node and need not be real. 
The EW structure $([g^\CC],\nabla)$ is induced so that, for each point $q\in S^+$, 
$\tilde p_2(\tilde p_1^{-1}(q))$ is a complex null surface.  
Notice that both $\tilde p_1$ and $\tilde p_2$ are holomorphic. 

The real EW manifold $W_0$ is contained in $W^\CC_0$ as a real slice. 
Let $i: W_0 \to W^\CC_0$ be the inclusion map. 
There is a natural identification map
$$ \mu : T^\CC W_0 \to TW^\CC_0|_{W_0} $$
defined by $\mu(v+ \sqrt{-1} w) = v+ J(w)$ where $J$ is the complex structure on $TW^\CC_0$
and $v, w$ are real tangent vectors on $W_0$. 
Recall that there is a natural identification $j : TW^\CC_0 \to T^{1,0}W^\CC_0$ 
defined by 
$ j(\xi) = \frac{1}{2}(\xi-\sqrt{-1}J(\xi))$
for $\xi \in  TW^\CC_0$. 
The real Einstein-Weyl structure $([g],\nabla)$ on $W_0$ is obtained as the restriction, that is the pullback by 
the identification  $j \circ \mu : T^\CC W_0 \to T^{1,0}W^\CC_0|_{W_0}$. 
Hence a complex null plane at $x\in W_0$ is considered as a complex null plane at $i(x)\in W^\CC_0$, 
so we obtain natural inclusion map $\iota : \mathcal N^+ \to \mathcal N^\CC$.

We introduce a holomorphic rank 2 distribution $\mathcal E \subset T^{1,0}\mathcal N^\CC$ so that, 
on $x=[v]\in \mathcal N^\CC$, $\mathcal E_x$ is the horizontal lift of the complex null plane 
$$ \Pi_x = \{ w \in T^{1,0}W_0^\CC \mid g(w,v)=0\}. $$ 
Since null surfaces on $W_0^\CC$ are totally geodesic, $\mathcal E$ is integrable and its integral manifolds are 
natural lift of null surfaces. 
Hence $\mathcal N^\CC$ is holomorphically foliated by such integral manifolds, 
and $S^+$ is its leaf space. 
Hence we obtain $\mathcal E = \ker\{ (\tilde p_1)_*: T^{1,0}\mathcal N^\CC \to T^{1,0}S^+\}$ 
or $j^{-1}(\mathcal E) =\ker\{ (\tilde p_1)_*: T\mathcal N^\CC \to TS^+\}$. 
By the construction, the complex distribution $E$ on $\mathcal N^+$ is the pull back of $\mathcal E$, 
that is, $E = \mu^{-1}(j^{-1}(\mathcal E))$. 

The key point is that the two maps $\mu: T^\CC W_0 \to T W_0^\CC|_{W_0}$ and 
$\iota_*: T^\CC W_0 \to T^\CC W^\CC_0$ do not coinside.  
Let us write $J$ for the complex structure on $\mathcal N^\CC$ and $S^+$. 
For $x\in \mathcal N^+$ and $u = v+ \sqrt{-1}w \in E_x$, 
we obtain $\mu(u) = v + J(w) \in j^{-1}(\mathcal E_x) = \ker (\tilde p_1)_*. $
Since $(\tilde p_1)_*$ is $J$-equivariant, we obtain 
\begin{equation} (\tilde p_1)_*(w) =J(\tilde p_1)_*(v).
\end{equation} 

On the other hand, $\iota_*(u) = v + \sqrt{-1}w \in T^\CC W^\CC_0$ and 
$$\begin{aligned}
(p_1)_*(u) &= (\tilde p_1)_*( \iota_*(v + \sqrt{-1}w))
 = (\tilde p_1)_*(v) + \sqrt{-1}(\tilde p_1)_*(w) \\
 &=  (\tilde p_1)_*(v) + \sqrt{-1} J  (\tilde p_1)_*(v)
\end{aligned} $$
Thus $(p_1)_*(u) \in T^{0,1}_{p_1(x)}S^+$. 
Hence \eqref{eq:01t} holds as required. 
\end{proof}

\subsection{Variation of the Einstein-Weyl structure}\label{ss:var2}
Recall from the explanation right before Remark \ref{rmk:Joyce} that the projection from any real point $P$ on the line $\ol{e_0e_1}$ realizes the Segre surface $S$ as a double covering of the cone $C(\Lmd)\subset\CP_3$ over the conic $\Lmd$, for which we denoted by $\Pi:S\lras C(\Lmd)$.
This preserves the real structure but does not preserve the $\CC^*$-action.
The branch divisor $\Sigma$ of $\Pi$ is a quartic which is a smooth elliptic curve not passing through the vertex of the cone.
This is real and we readily see that the real locus $\Sigma\us$ consists of two circles.
The images $\mathbb D_1:=\Pi(S\us_1)$ and $\mathbb D_2:=\Pi(S\us_2)$ of the real spheres are closed disks in $C(\Lmd)$ that are mutually disjoint, and $\Pi|_{S\us_i}$ ($i=1,2)$ are double covering over these disks branched along the boundary circles.
These circles are exactly the real locus of $\Sigma$.
We denote $\Sigma\us_i$ for these. So $\ptl\mathbb D_i = \Sigma\us_i$ for $i=1,2$.

Let $C\in W_0$ be any regular minitwistor line and $H\subset\CP_4$ the real hyperplane such that $C = S\cap H$.
Since $C$ is regular, $e_0,e_1\not\in H$.
Using the real point $P=\ol{e_0e_1}\cap H$ as the center of the projection,
we obtain the double covering $\Pi:S\lras C(\Lmd)$ which satisfies  
$C=\Pi\inv(h\cap C(\Lmd))$ where $h=\Pi(H)$, a real hyperplane not passing through the vertex of the cone.
Because $C$ has exactly one node, the plane $h$ is tangent to the branch elliptic curve $\Sigma$ at exactly one point and the node of $C$ is over the tangent point.
Hence, this point belongs to the circle $\Sigma\us_1$.
The residual two points of $\Sigma\cap h$ belong to another circle $\Sigma\us_2$, and they are the endpoints of the intersection $h\cap \mathbb D_2$ which is a closed arc in the disk.
The locus $\Pi\inv(h\cap \mathbb D_2)$ is the real circle $C^{\RR}=C\cap S\us_2$ of $C$.

Next, take any irregular minitwistor line.
It is of the form $l_i + \ol l_i + D_{\lmd}$, where $i=1,2$, $\lmd\in I_2=(\lmd_3,\lmd_4)\subset\Lmd\us$ and $D_{\lmd}\subset S$ is the conic over $\lmd$.

\begin{proposition}\label{p:modif1}
If $C$ and $l_i + \ol l_i + D_{\lmd}$ are respectively any regular and irregular minitwistor lines on $S$ as above, then there exists a real open set $U$ in $S$ containing both $C$ and $l_i + \ol l_i + D_{\lmd}$, such that $U$ admits the modification given in Section \ref{ss:modify}.
\end{proposition}

\proof For simplicity,
we denote $C_0:= l_i + \ol l_i + D_{\lmd}$ and use the above notations $\Pi$, $\Sigma$, $H$ and $h$, so that $C = S\cap H = \Pi\inv(C(\Lmd)\cap h)$.
Further, we identify any branch point of $\Pi$ with the ramification point over it.
Let $p\in S\us_1$ be the node of $C$.
Let $h_0\subset\CP_3$ be the hyperplane that satisfies $C_0=\Pi\inv(C(\Lmd)\cap h_0)$. This passes through the vertex of $C(\Lmd)$. 
Let $\delta$ be the intersection point of the generating line of $C(\Lmd)$ over $\lmd$ and the plane $h$.
This is a real point and belongs to the interior of the disk
$\mathbb D_2$, so 
$\Pi\inv(\delta)$ consists of two points of $S\us_2$.
Let $x$ and $y$ be these points.
For any point $\zeta\in \Sigma\us_1$, the tangent line $T_{\zeta}\Sigma$ and the point $\delta$ span a hyperplane.
We denote $h(\zeta;\ddd)\subset\CP_3$ for it
and denote $C(\zeta;\ddd):=\Pi\inv(h(\zeta;\ddd)\cap C(\Lmd))$, which is a minitwistor line that has a node at $\zeta$ and which passes through $x$ and $y$.
Note that this definition of $C(\zeta;\ddd)$ makes sense even when $\zeta\in\Sigma$ is non-real as a curve in $S$ as long as $T_{\zeta}\Sigma$ and $\ddd$ span a hyperplane.
If $\mf C_{x,y}$ is the spacelike geodesic formed by minitwistor lines through $x$ and $y$ as in Section \ref{ss:slg},
then 
$$
\mf C_{x,y} = \big\{C(\zeta;\ddd)\set \zeta\in\Sigma\us_1\big\}.
$$
Hence, there is an identification $\Sigma\us_1\simeq\mf C_{x,y}$ which assigns $C(\zeta;\ddd)$ to each $\zeta\in\Sigma\us_1$.
Under this identification, the points $C_0, C\in \mf C_{x,y}$ correspond to the points $p_i$ and $p$ respectively.

In the following we let $i=1$ for simplicity of presentation and
let $A(p_1, p)$ be the closed arc in $\Sigma\us_1$ whose ends are $p_1$ and $p$ and which does not contain $p_2$.
Then for any distinct $\zeta,\zeta'\in A(p_1, p)$, 
one of the two intersection points of the plane sections $h(\zeta;\ddd)\cap C(\Lmd)$ and $h(\zeta';\ddd)\cap C(\Lmd)$ belong to $f\inv(\ol I_1)$, where $I_1$ is the arc $(\lmd_1,\lmd_2)$ in the circle $\Lmd\us$ as before, $\ol I_1$ is its closure, and $f:C(\Lmd)\lras \Lmd$ is the cone projection.
Therefore, for any open tubular neighborhood $T\subset\Lmd$ of the closed arc $f(A(p_1,p))$,
there exists an open tubular neighborhood $T'$ in $\Sigma$ of $A(p_1, p)$, such that 
\begin{align}\label{triple}
\big(h(\zeta;\ddd)\cap h(\zeta';\ddd)\cap C(\Lmd)\big)\minus\{\ddd\}\in f\inv(T) {\text{ for any $\zeta\neq\zeta'\in T'$.}}
\end{align}
Note that the triple intersection consists of two points as $\deg C(\Lmd)=2$.
This intersection is of course transverse as an intersection of two hyperplane sections.
We choose and fix the open neighborhood $T$ of $f(A(p_1,p))$ in such a way that it does not intersect a neighborhood of the closed arc $\ol I_2=[\lmd_3,\lmd_4]$ in $\Lmd$.
In the following we take both $T$ and $T'$ to be real.

As in Definition \ref{d:pm}, as long as $\zeta\in \Sigma\us_1$, we have the decomposition $C(\zeta;\ddd)=C(\zeta;\ddd)^+\cup C(\zeta;\ddd)^-$ into disks, where $C(\zeta;\ddd)^+$ is the part that intersects the line $\ol l_1$ (and also $l_2$).
Shrinking the neighborhood $T'$ of $A(p_1,p)$ if necessary, 
we can smoothly extend this decomposition for any $\zeta$ belonging to $T'$ in such a way that $x,y\in C(\zeta;\ddd)^+\cap C(\zeta;\ddd)^-$ and they are compatible with the real structure in the sense $\sigma(C(\zeta;\ddd)^+) = C(\sigma(\zeta);\ddd)^-$.

Next, we take a tubular neighborhoods $L_1$ and $L_2$ of $\ol I_2$ in $\Lmd$ such that $\ol L_1\subset L_2$, $L_2\cap T=\emptyset$ and $L_1, L_2$ real. Then $x,y\in f\inv(L_1)\subset f\inv(L_2)$.
Using $L_1$ and the above neighborhood $T'\supset A(p_1,p)$ in $\Sigma$, we define two subsets $U^+$ and $U^-$ of $S$ by 
\begin{align}\label{}
U^{+}=\Big(\bigcup_{\zeta\in T'} C(\zeta;\delta)^+\Big)\setminus f\inv(\ol L_1)
\qandq
U^{-}=\Big(\bigcup_{\zeta\in T'} C(\zeta;\delta)^-\Big)\setminus  f\inv(\ol L_1).
\end{align}
Obviously, $\sigma(U^+) = U^-$.
From \eqref{triple}, noting that the left-hand side of \eqref{triple} necessarily belongs to $C(\zeta;\ddd)^+\cap C(\zeta;\ddd)^-$ under the identification of ramification points and branch points, for any $\zeta\neq\zeta'\in T'$,
we have
\begin{align}\label{transv}
\big(C(\zeta;\delta)^+ \cap C(\zeta';\delta)^-\big)\minus\{x,y\} \subset f\inv(T).
\end{align}
Therefore, noting that the real curve $(C(\zeta;\delta)^+ \cap C(\zeta;\delta)^-)\minus\{\zeta\}$ stays in a neighborhood of the second sphere $S\us_2$ in $S$ if $\zeta\in T'$ and therefore it is disjoint from $f\inv(T)$ as $L_2\cap T=\emptyset$,
\begin{align}\label{pmT}
U^+\cap U^-\subset f\inv(T).
\end{align}
Moreover, as above, the left-hand side of \eqref{transv}, which consists of two points over $\big(h(\zeta;\ddd)\cap h(\zeta';\ddd)\cap C(\Lmd)\big)\minus\{\ddd\}$, are transversal intersection.
Hence, by taking a thin tubular neighborhood of each $C(\zeta;\delta)$ and taking a `normalization' of it as in the procedure from $\tilde U$ to $\hat U$ in Section \ref{ss:modify}, the transformations of the curves $C(\zeta;\ddd)$ into the `normalization' intersect only at $x$ and $y$ as $C(\zeta;\delta)^2 = \deg S = 4$. 
This implies that $U^+$ and $U^-$ are open subsets of $S$.

From \eqref{pmT}, gluing the disjoint union $U^+\sqcup U^-$ and $f\inv(L_2)$ by identifying the open subset $(U^+\sqcup U^-)\cap f\inv(L_2)
= (U^+\cap f\inv(L_2)) \sqcup(U^-\cap f\inv(L_2))$ of the two sets,
we may define 
\begin{align}\label{hU}
\hat U:= (U^+\sqcup U^-) \cup f\inv(L_2)
\end{align}
which is real and has a natural projection to the open subset $U:=(U^+\cup U^-) \cup f\inv(L_2)$ of $S$.
Then $U$ contains $C_0$ and $C$, and the minimal resolution of \eqref{hU} at $e_0$ and $e_1$ clearly plays the role of $\hat U$ in the construction of Section \ref{ss:modify}.
\proofend

\begin{figure}
\includegraphics{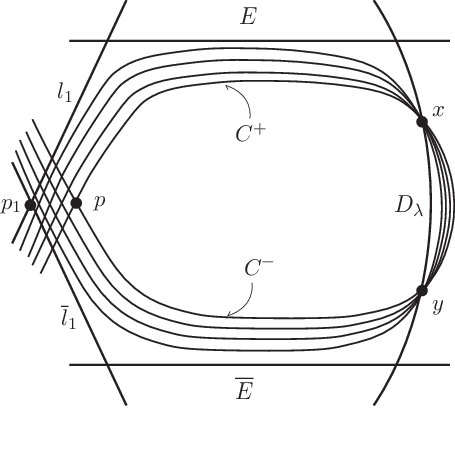}
\caption{
The family $\{C(\zeta;\ddd)\}_{\zeta\in T'}$ in the resolution 
}
\label{fig:modif2}
\end{figure}

\medskip
For the family $\{C(\zeta;\ddd)\}_{\zeta\in T'}$ in the above proof, see Figure \ref{fig:modif2}, where the minimal resolution of $S$ is taken.

Using Proposition \ref{p:modif1}, we show the next proposition, which will be soon used.
\begin{proposition}\label{p:aut1}
Let $S$ and $S'$ be Segre surfaces defined by the equation \eqref{Segre1} using the real numbers $0<\kappa,\kappa'<1$ respectively, and $W$ and $W'$ the EW spaces obtained from $S$ and $S'$ respectively.
If $\phi:W\lras W'$ is any isomorphism as EW spaces, then 
$\phi(W_0) = W'_0$.
\end{proposition}

\proof  Let $C\in W_0$ be any regular minitwistor line on $S$ and $l_i + \ol l_i + D_{\lmd}$ ($i=1,2$ and $\lmd\in I_2$) be any irregular minitwistor line on $S$. From Proposition \ref{p:modif1}, there exists
an open subset $U\subset S$ containing both $C$ and $l_i + \ol l_i + D_{\lmd}$
for which the modification can be applied, 
and let $\widecheck U$ be the modification of $U$. 
Let $\check C$ be the transformations of $C$ to $\check U$. By the same argument as in Proposition \ref{p:modmtl}, by making use of the transformation of $l_i + \ol l_i + D_{\lmd}$ into $\check U$, we obtain that 
$\check C$ is a smooth minitwistor line in $\check U$ and it satisfies $\dim |\check C| = 0$.

On the other hand, let $l'_i + \ol l'_i + D'_{\lmd'}$, $i=1,2$ and $\lmd'\in I'_2=(\lmd'_3,\lmd'_4)$, be any irregular minitwistor line on $S'$ and $U'\subset S'$ an open subset that contains $l'_i + \ol l'_i + D'_{\lmd'}$ which is chosen as we did in Section \ref{ss:modify}.
Let $\check U'$ be the modification of $U'$ and $\check D'_{\lmd}$ the tranformation of $l'_i + \ol l'_i + D'_{\lmd'}$ into $\check U'$.
By Proposition \ref{p:modmtl}, $\dim |\check D'_{\lmd'}| = 1$.

Since the EW structure on $W$ (resp.\,$W'$) around the point $C\in W_0$ (resp.\, $l'_i + \ol l'_i + D'_{\lmd'}\in W'\minus W'_0$) is determined from a neighborhood of $\check C$ (resp.\,$\check D'_{\lmd'}$) in $\check U$ (resp.\,in $\check U'$), this implies that the EW structure around the point $C$ can never be isomorphic to that around the point $l'_i + \ol l'_i + D'_{\lmd'}\in W'\minus W'_0$.
Therefore, $\phi(C)\in W'_0$. Hence, $\phi(W_0) = W'_0$.
\proofend

\medskip
Using this as well as the result in the previous subsection, we show the following, which is one of the main results in this section.
\begin{theorem}\label{t:var}
Let $S,S'$ and $W,W'$ be as in Proposition \ref{p:aut1}.
Then $W\simeq W'$ as EW spaces if and only if there exists a holomorphic isomorphism $S\simeq S'$ which commutes with the real structure.
\end{theorem}

\proof
The implication $S\simeq S'\Longrightarrow W\simeq W'$ is obvious. Let $\phi:W\lras W'$ be an isomorphism as EW spaces. 
By Proposition \ref{p:aut1}, $\phi$ maps $W_0$ to $W'_0$. 
Then from Corollary \ref{c:01t}, this induces a holomorphic automorphism $\phi^+:S^+\lras S^+$.

On the minimal resolution $\tilde S$ of $S$, $l_1$ and $\ol l_2$ are $(-1)$-curves.
Let $\tilde S\lras \tilde S_1$ be the blowdown of these two curves.
Then the images of $E$ and $\ol E$ are $(-1)$-curves on $\tilde S_1$.
Let $\tilde S_1\lras \tilde S_2$ be the blowdown of these two curves.
It is easily seen that $\tilde S_2$ is biholomorphic to $\qdr$.
Let $e$ and $\ol e\in \tilde S_2$ be the images of $E$ and $\ol E$.
These points belong to neither the same $(1,0)$-curve nor the same $(0,1)$-curve.
Similarly, from another surface $S'$, we obtain the analogous sequence $\tilde S' \lras \tilde S'_1\lras\tilde S'_2\simeq\qdr$ and two points $e'$ and $\ol e'$ on $\tilde S'_2$.
These contractions do not preserve the real structure.

Now from the above construction, for the surface $S^+=S\minus( l_1\cup \ol l_2)$ and $(S')^+=S'\minus(l'_1\cup \ol l'_2)$, we have obvious holomorphic identifications
\begin{align}\label{har1}
S^+=S\minus (l_1\cup \ol l_2)\simeq 
\tilde S\minus(l_1\cup \ol l_2\cup E\cup \ol E)
\simeq \tilde S_1\minus (E\cup \ol E)
\simeq \tilde S_2\minus \{e,\ol e\}
\end{align}
and similarly for $(S')^+$.
Therefore, the biholomorphism $\phi^+:S^+\lras (S')^+$
induces a biholomorphism
$\phi_2^+:\tilde S_2\minus \{e,\ol e\}\lras \tilde S'_2\minus \{e',\ol e'\}$.
By Hartogs' theorem, $\phi_2^+$ extends to a biholomorphism $\phi_2:\tilde S_2\lras\tilde S'_2$ which maps $\{e,\ol e\}$ to $\{e',\ol e'\}$.
This means that $\phi_2$ maps $(1,1)$-curves through $e,\ol e$ to $(1,1)$-curves through $e',\ol e'$.
Hence, the lift $\phi_1:\tilde S_1\lras\tilde S'_1$ maps fibers of the projection to $\Lmd$ to fibers of the projection to $\Lmd'$, and maps the fibers $\ol l_1\cup l_2$ to the fibers $\ol l'_1\cup l'_2$. This implies that $\phi_1$ maps the image points of $l_1$ and $\ol l_2$ under the blowing down $\tilde S\lras\tilde S_1$ to the image points of $l'_1$ and $\ol l'_2$ under the blowing down $\tilde S'\lras\tilde S'_1$
Hence, $\phi_1$ lifts to an isomorphism $\tilde\phi:\tilde S\lras \tilde S'$ and it satisfies $\tilde\phi(E\cup\ol E) = E'\cup\ol E'$.  
Therefore, $\tilde\phi$ descends to a holomorphic isomorphism $\phi:S\lras S'
$, which gives the extension of $\phi^+$.

It remains to show that $\phi$ preserves the real structure.
The real structure on $S$ is defined by thinking of each point on it as a complex null direction $[v]$ at a point of $W_0$ and then assigns the complex conjugation $[\ol v]$ at the same point of $W_0$.
Since the differential of $\phi$ is of course defined over $\RR$, this readily means the required commutativity regardless of whether the differential preserves $S^+$ or flips $S^+$ and $S^-$.
\proofend

\medskip
The complex structure of the Segre surface $S$ depends on 
the cross-ratio of the four points $\lmd_1,\lmd_2,\lmd_3,\lmd_4$ on $\Lmd\simeq\CP_1$ and it is $(1-\kappa)\inv$.
Since $0<\kappa<1$, from Theorem \ref{t:var}, this implies that the moduli space of the EW structure on $W\simeq\mathbb S^2\times I$ is identified with the interval $(0,1)$.
Further, the proof of the theorem implies that the group of automorphism of the EW structure on $W$ associated to a Segre surface $S$ is naturally isomorphic to the group of holomorphic automorphism group of $S$ preserving the real structure. 
Since the identity component of the latter automorphism group is the subgroup $S^1\subset \CC^*$ as remarked in Section \ref{s:S}, we obtain the following 
\begin{theorem}\label{t:aut}
The identity component of the automorphism group of the
Lorenzian Einstein-Weyl structure on $W$ is a circle.
\end{theorem}

As remarked in Section \ref{s:I}, since the implication ``EW space $\Longrightarrow$ minitwistor space'' is established only locally, Theorems \ref{t:var} and \ref{t:aut} cannot be immediately obtained from the relevant properties of the Segre surfaces.

\subsection{Degeneration to de Sitter space}\label{ss:dS}
As above, the Segre surface $S\subset\CP_4$ in \eqref{Segre1} depends on the parameter $\kappa\in (0,1)$. In this subsection, we prove that in the limit $\kappa\to 0$, the EW structure on $W$ converges to the de Sitter structure. 
This proves the second assertion in Theorem \ref{t:2}.

For this, we first write down the equation of arbitrary real hyperplane $H\in W$ explicitly.
Using spherical coordinates $(\rho,\theta)\in[-1.1]\times [0,2\pi)$,
from the description of the sphere $S\us_1$ in the proof of Proposition \ref{prop:Ssigma}, we may represent 
points on $S\us_1$ by
$$
(X_0,X_1,X_2,X_3,X_4) = \left(\sqrt{1-\rho^2}e^{i\theta},\sqrt{1-\rho^2}e^{-i\theta}, \frac 1{\sqrt{\kappa}},\, \rho,\,\xi(\rho)\right) =:p(\rho,\theta),
$$
where $\xi(\rho) = \sqrt{\frac 1\kappa-\rho^2}$.
Here, $\theta$ parametrizes $\mathbb S^1$-orbits on $S\us_1$ and the fixed points $p_1,p_2$ correspond to the cases $\rho=\pm 1$.
Using this, we obtain that real hyperplanes in $\CP_4$ which are tangent to $S$ at the point $p(\theta,\rho)\in S\us_1$ are defined by 
\begin{align}\label{explct_hypsf}
\frac{\sqrt{1-\rho^2}}2\big(e^{-i\theta}X_0 + e^{i\theta}X_1\big)
- \sqrt{\kappa}X_2 + \rho X_3
+ \mu\left(
\frac 1{\sqrt\kappa} X_2 - \rho X_3 - \xi X_4
\right)
=0,
\end{align}
where $\mu\in\RR\cup\{\infty\}= (T_{p(\theta,\rho)}S)^{*\sigma}$.
We write the hyperplane \eqref{explct_hypsf} by $H(p,\mu)$.
The exceptional hyperplanes $H_2,H_3$ and $H_4$ given in Definition \ref{d:ehpls} occur when the corresponding coefficient vanishes, that is when $\mu=\kappa, 1$ and $0$ respectively.
By the previous argument in Section \ref{s:rmtl} (see Figure \ref{fig:circle2}), 
the interval $J_4(p)\subset (T_{p(\theta,\rho)}S)^{*\sigma}$ is parameterized by $\mu\in(0,\kappa)$.
The triple 
\begin{align}\label{Wcoord}
(\rho,\theta,\mu)\in [-1,1]\times[0,2\pi)\times (0,\kappa)
\end{align}
can be used as global coordinates on $W\simeq\mathbb S^2\times I$.

Next, using the parameter $\kappa\in(0,1)$, we shall introduce new homogeneous coordinates on $\CP_4$ by
\begin{align}\label{change_X_Z}
(Z_0,Z_1,Z_2,Z_3,Z_4) = (X_0,X_1,\sqrt\kappa X_2,X_3,X_2 + X_4).
\end{align}
Then the real structure $\sigma$ and the $\CC^*$-action take the same forms as those for the coordinates $(X_i)$
(see \eqref{rs1} and \eqref{C*action1}), while the equation of $S=S(\kappa)$ is written as 
\begin{align}\label{Segre2}
Z_0Z_1 - Z_2^2 + Z_3^2 = \sqrt\kappa(Z_3^2+Z_4^2) - 2Z_2Z_4=0.
\end{align}
We consider the limit $\kappa\searrow 0$. Then the equations converge to 
\begin{align}\label{def_T}
Z_0Z_1 - Z_2^2 + Z_3^2 = Z_2Z_4=0.
\end{align}
Let $T\subset\CP_4$ be the surface defined by \eqref{def_T}.
This consists of two irreducible components $T_1\subset\{Z_2=0\}$ and $T_2\subset\{Z_4=0\}$.
$T_1$ is a cone over a conic with the vertex being $p_{\infty}:=(0,0,0,0,1)$ and $T_2$ is a smooth quadric.
Further, $T_1$ and $T_2$ are $\sigma$-invariant, $\mathbb S^1$-invariant, $T_1\us = \{p_{\infty}\}$ and $T_2\us\simeq \mathbb S^2$.
In particular, $T_2$ is the minitwistor space of the de Sitter space, where their minitwistor lines are hyperplane sections that intersect $T\us_2$ in a circle. Further, $T\us_2$ is $\mathbb S^1$-invariant.

Recall that the $\CC^*$-fixed locus on the Segre surface $S=S(\kappa)$ consist of the 6 points $e_0,e_1,p_1,p_2,p_3,p_4$. 
The singular points $e_0$ and $e_1$ are independent on $\kappa$ and belong to $T_1\cap T_2$, which is a smooth conic.
Both the points $p_1,p_2\in S\us_1$ converge to the vertex $p_{\infty}\in T_1$ in the limit $\kappa\to 0$.
In contrast, another pair of $\mathbb S^1$-fixed points $p_3,p_4\in S\us_2$ converge to different points of $T\us_2$:
$$
\lim_{\kappa\to 0} p_3 = (0,0,1,-1,0)
\qandq
\lim_{\kappa\to 0} p_4 = (0,0,1,1,0).
$$
We write $p'_3$ and $p'_4$ for these points in the limit.
The $\CC^*$-fixed points on $T$ consist of the 5 points
$e_0,e_1,p'_3,p'_4$, and $p_{\infty}$.

The projection $f:\CP_4\lras\CP_2$ from the line $\ol{e_0e_1}$ is still the map that drops the first two coorinates $Z_0$ and $Z_1$, and from \eqref{Segre2} 
the conic $\Lmd=f(S)$ is written as
$$
\Lmd = \left\{
(Z_2,Z_3,Z_4)\set  \sqrt\kappa(Z_3^2+Z_4^2) - 2Z_2Z_4=0
\right\}.
$$
Recall that $\Lmd$ contains the real four points $\lmd_i = f(p_i)$
$(i=1,2,3,4)$ as the discriminant points of the conic bundle. 
Taking the limit $\kappa\to 0$, we define
$$
\Theta = \left\{
(Z_2,Z_3,Z_4)\set  Z_2Z_4=0
\right\},
$$
which consists of two lines $\Theta_1 = \{Z_2=0\}$ and  $\Theta_2 = \{Z_4=0\}$. Then $f(T_i) = \Theta_i$ for $i=1,2$.
The $\CC^*$-action \eqref{C*action1} survives in the limit and $f|_{T_i}:T_i\to\Theta_i$ ($i=1,2$) can be regarded as a quotient map for the $\CC^*$-action on the limit.
Let 
$
\lmd_0:= (0,1,0)
$
be the intersection point $\Theta_1\cap \Theta_2$.
Then $f(T_1\cap T_2) = \{\lmd_0\}$.
The points $\lmd_3,\lmd_4\in \Lmd$ respectively converge to 
$$
\lmd_3' = (1,-1,0) = f(p_3)
\qandq
\lmd_4' = (1,1,0) = f(p_4),
$$
while the points $\lmd_1,\lmd_2\in\Lmd$ converge to the same point 
$$
\lmd_{\infty} = (0,0,1) = f(p_{\infty}).
$$
The arc $I_1\subset\Lmd\us$ bounded by $\lmd_1$ and $\lmd_2$ shrinks to the point $\lmd_{\infty}$ and another arc $I_2\subset\Lmd\us$ bounded by $\lmd_3$ and $\lmd_4$ converges to the arc in $\Theta_2\us$ bounded by $\lmd'_3$ and $\lmd'_4$.
It follows that, {\em in the limit $\kappa\to 0$, the sphere $S\us_1$ shrinks to the vertex $p_{\infty}$ of the cone $T_1$, while another sphere $S\us_2$ converges to the whole sphere $T\us_2$
in the minitwistor space $T_2$ of the de Sitter space.}
In a similar way, we also obtain that, in the limit $\kappa \to 0$,
the two lines $l_1,l_2\subset S$ converge the common generating line of the cone $T_1$,
and the two lines $\ol l_1,\ol l_2\subset S$ converge another generating line in $T_1$.
Let $l_{\infty}$ and $\ol l_{\infty}$ be the former and the latter generating lines in $T_1$ respectively.
Further, let $l_3'$ and $l_4'$ be the limit lines of $l_3$ and $l_4$ as $\kappa\to 0$. These are lines in $T_2$ and components of the fibers of $f$ over the limit points $\lmd_3'$ and $\lmd_4'\in \Theta_2$ respectively.

We now define an involution $\aaa$ on the space 
$[-1,1]\times[0,2\pi)\times I$ with $I=(0,\kappa)$ by 
\begin{align}\label{aaa}
(\rho,\theta,\mu) \stackrel{\aaa}\longmapsto 
(-\rho,\theta+\pi,\kappa - \mu),
\end{align}
where the middle component is considered modulo $2\pi$.
This has no fixed point and can be thought of as an involution on the space $W\simeq\mathbb S^2\times I$ because $\aaa$ maps the hyperplane $H(p,\mu)$ to $H(q,\kappa-\mu)$, where $p=(\rho,\theta)$ and $q = (-\rho,\theta+\pi)$. 
Here, we are not asserting that $\aaa$ preserves the EW structures on $W$.

By putting $\mu=\kappa\DDD$ and taking the limit $\kappa\to0$ in \eqref{explct_hypsf}, 
we obtain that the hyperplane $H(p,\mu)$ 
converges to the hyperplane in $\CP_3$ defined by 
\begin{align}\label{limhp}
\frac{\sqrt{1-\rho^2}}2\big(e^{-i\theta}Z_0 + e^{i\theta}Z_1\big)
+\rho Z_3
=(1-2\DDD)Z_2, \quad 0<\DDD<1.
\end{align}
Note that if we replace $(\rho,\theta,\DDD)$ with $\aaa(\rho,\theta,\DDD)$ in this equation, then both sides will be multiplied by $(-1)$, which means that 
$\lim_{\kappa\to 0} H = \lim_{\kappa\to0} \aaa(H)$ and that the triple $(\rho,\theta,\DDD)$ do not effectively parametrize these limit hyperplanes.
Further, by elementary calculation, the coincidence in the limit happens only for such cases.

Because the sphere $S\us_2$ converges to the whole sphere $T\us_2$ as above and any hyperplane of the form \eqref{explct_hypsf} intersects $S\us_2$ in a circle, the limiting hyperplane \eqref{limhp} always intersects the sphere $T\us_2$ in a circle.
Let $W'_{\rm dS}$ be the space of real hyperplanes in $\CP_3$ that intersect  $T_2\us$ in a circle.
It can also be seen by elementary calculation that the hyperplanes of the form \eqref{limhp} (where $|\rho|\le 1$, $\theta\in[0,2\pi)$ and $\DDD\in (0,1)$) are all such real hyperplanes in $\CP_3$.
For any $H\in W'_{\rm dS}$, the section $T_2\cap H\simeq\CP_1$ is divided into disks by $T\us_2$, and if we denote $W_{\rm dS}$ for the space of the disks obtained this way, then $W_{\rm dS}$ is the de Sitter space and there is an unramified double covering map $W_{\rm dS}\lras W'_{\rm dS}$. $W_{\rm dS}$ is diffeomorphic to $\mathbb S^2\times I$ and this map can be regarded as a universal covering map for $W'_{\rm dS}$.

Thus, there are twice as many minitwistor lines in $S$ as those in the quadric $T_2$, and the same for the spaces of disks.
The next proposition means that the $+$-disks in the hyperplane sections $S\cap H$ and $S\cap \aaa(H)$ converge to {\em distinct} disks to form one minitwistor line in $T_2$.
(The same thing holds for the $-$-disks.)

\begin{proposition}
For any $H\in W$, let $C=S\cap H$ and $C=C^+\cup C^-$ be the decomposition into
$+$-part and $-$-part as discussed in Section \ref{ss:disk}, so that $C^+\subset \ms C^+$ and $C^-\subset \ms C^-$.
Similarly, for $\aaa(H)\in W$, let $\aaa(C) = \aaa(H)\cap S$, and 
$\aaa(C)^+$ and $\aaa(C)^-$ the $+$-part and $-$-part of $\aaa(C)$ respectively.
Then in the limit $\kappa\to 0$, 
$C^+$ and $\aaa(C)^+$ converge to distinct disks in the limit minitwistor line in $T_2$. Namely,
$$
\left(\lim_{\kappa\to 0} C^+\right)
\cup
\left(\lim_{\kappa\to 0} \aaa(C)^+\right)
= 
\left(\lim_{\kappa\to 0} H\right)\cap T_2.
$$
\end{proposition}

\proof 
Evidently, the property $\lim_{\kappa\to 0}C^+\neq \lim_{\kappa\to 0}\aaa(C)^+$ is open and closed condition. Therefore, since $W$ is connected, to prove that this property holds for any $C=S\cap H\in W$, it suffices to show for some $C=S\cap H\in \ol W$.
We choose a very specific member $H=H_2(p)\in \ms H_2$.
Then $S\cap H_2(p_1) = l_1+\ol l_1 + l_3 + \ol l_3$.
In this curve, the half $l_1+l_3$ is the $+$-part since this intersects $\ol l_1$ at a point, see Definition \ref{d:pm}.
From the definition $\aaa$ as in \eqref{aaa}, 
we readily obtain $\aaa(H_2(p_1)) = H_4(p_2)\in\ms H_4$.
Then $S\cap H_4(p_2) = l_2+\ol l_2 + l_3 + \ol l_3$.
The $+$-part of this curve is $\ol l_2 + \ol l_3$ since this intersects $l_2$ at a point.

If we take the limit $\kappa\to 0$, then the $+$-part $l_1+l_3$ of $S\cap H_2(p_1)$ converges to $l_{\infty} + l_3'$, while the $+$-part $\ol l_2+ \ol l_3$ of $S\cap H_4(p_2)$ converges to $\ol l_{\infty} + \ol l_3'$.
Thus, the limits in $T_2$ are $l'_3$ and $\ol l'_3$, and they are indeed distinct, divided by the real point $l'_3\cap \ol l'_3\in T\us_2$.
\proofend

\medskip
As an immediate consequence of the proposition, we obtain:

\begin{corollary}
In the limit $\kappa\to 0$, the EW structure on $W$ obtained from $S$ converges to the de Sitter structure.
\end{corollary}

This proves the second assertion in Theorem \ref{t:2}.


\medskip

\end{document}